\documentclass[11pt,oneside, final]{amsart}
\copyrightinfo{0}{Iranian Mathematical Society}
\pagespan{1}{\pageref*{LastPage}}
\usepackage{etoolbox,lastpage}
\commby{}
\usepackage{amsmath,amsthm,amscd,amsfonts,amssymb,enumerate}
\usepackage{graphicx}
\usepackage{color}
\usepackage[colorlinks]{hyperref}

\makeatletter \@addtoreset{equation}{section}

\newtheorem{theorem}{Theorem}[section]
\newtheorem{proposition}[theorem]{Proposition}
\newtheorem{lemma}[theorem]{Lemma}
\newtheorem{corollary}[theorem]{Corollary}
\newtheorem{conjecture}[theorem]{Conjecture}
\theoremstyle{definition}
\newtheorem{definition}[theorem]{Definition}

\newtheorem{remark}[theorem]{Remark}
\newtheorem{assumption}[theorem]{Assumption}
\theoremstyle{remark}
\numberwithin{equation}{section}

\begin{document}


\title[Mean curvature flow \& Micheal-Simon
inequalities]{Mean curvature type flow and sharp Micheal-Simon
inequalities}

\author[J. Cui]{Jingshi Cui}

\author[P. Zhao]{Peibiao Zhao$^*$}

  \thanks{$^*$Corresponding author}
%
\begin{abstract}
In this paper, we first investigate a new locally constrained mean
curvature flow (\ref{1.5}) and prove that if the initial
hypersurface $M_{0}$ is of smoothly compact starshaped, then the
solution $M_{t}$  of the flow (\ref{1.5}) exists for all time and
converges to a sphere in $C^{\infty}$-topology. Following this flow
argument, not only do we achieve a new proof of the celebrated sharp
Michael-Simon inequality for mean curvature in Euclidean space
$\mathbb{R}^{n+1}$, but we also get the necessary and sufficient
condition for the establishment of the equality.

In the second part of this paper, we study a mean curvature type
flow (\ref{1.7}) of static convex hypersurfaces in Euclidean space
$\mathbb{R}^{n+1}$, and prove that the flow (\ref{1.7}) has a unique
smooth solution $M_{t}$ for all time $t\in[0,+\infty)$, and the
static convexity of the hypersurface is preserved along the flow
(\ref{1.7}). Moreover, $M_{t}$ converges exponentially to a sphere
of radius $R$ in $C^{\infty}$-topology as $t \to +\infty$. By
exploiting the properties of this flow, we develop and present a new
sharp Michael-Simon inequality for $k$th mean curvature.\\
\textbf{Keywords:}   Locally constrained mean curvature flow; Mean curvature type flow; Michael-Simon inequality; $k$th mean curvature.  \\
\textbf{MSC(2010):}  Primary: 53E99; Secondary: 52A20; 35K96
\end{abstract}

 \maketitle
%

\section{\bf Introduction}

\noindent It is well known that the classical Michael-Simon
inequality as a Sobolev inequality on submanifolds of Euclidean
space  is not sharp. The classical Michael-Simon inequality says as
follows
\begin{theorem}{(\cite{MS1973})}\label{th1.1}
    Let $i: M^{n} \rightarrow \mathbb{R}^{N}$ be an isometric immersion $(N>n)$. Let $U$ be an open subset of $M^n$. For a function $\phi  \in C_{c}^{\infty}(U)$, there exists a constant $C$, such that
    \begin{align}\label{1.1}
        \left(\int_{M}|\phi|^{\frac{n}{n-1}} d \mu_{M}\right)^{\frac{n-1}{n}} \leq C \int_{M}(|H| \cdot|\phi|+|\nabla \phi|) d \mu_{M}
    \end{align}
     where $H$ is the mean curvature of $M^n$.
\end{theorem}
If $\phi \equiv 1$, the Michael-Simon inequality gives the
relationship between the area and the integral of  mean curvature.
The best constant problem  in (\ref{1.1}) is still  an open issue,
even in the case of minimal surfaces. In a very recent paper,
Brendle \cite{B} had confirmed   a sharp version of the
Michael-Simon inequality as below.
\begin{theorem}{(\cite{B})}\label{th1.2}
    Let $M$ be a compact hypersurface in $\mathbb{R}^{n+1}$ (possibly with boundary $\partial M$ ), and let $f$ be a positive smooth function on $M$. Then
    \begin{align}\label{1.2}
       \int_{M} \sqrt{\left|\nabla^{M} f\right|^{2}+f^{2} H^{2}}+\int_{\partial M} f \geq n\left|B^{n}\right|^{\frac{1}{n}}\left(\int_{M} f^{\frac{n}{n-1}}\right)^{\frac{n-1}{n}}
    \end{align}
    where $H$ is the mean curvature of $M$ and $|B^{n}|$ is the area of the unit sphere $\mathbb{S}^{n}$ in $\mathbb{R}^{n+1}$. Moreover, if equality holds, then $f$ is constant and $M$ is a flat disk.
\end{theorem}
The method of proving the inequality (\ref{1.2}) in \cite{B} is
inspired by the Alexandrov-Bakelman-Pucci maximum principle (see
e.g. \cite{XC},  \cite{NT}) and the key step in the proof is to
simplify the original inequality (\ref{1.2}), which can be divided
into the following three folds.

The first fold is to derive
 that there holds
\begin{align*}
    \int_M \sqrt{|\nabla^{M} f|^{2}+f^{2}H^{2}}d\mu + \int_{\partial M}f d\mu = n\int_{M} f^{\frac{n}{n-1}}
    d\mu;
\end{align*}

The second fold is to prove that there is a function $\psi : M
\rightarrow \mathbb{R}$ which solves the PDE
\begin{align*}
    \operatorname{div}_{M}\left(f \nabla^{M} \psi \right)=n f^{\frac{n}{n-1}}-\sqrt{\left|\nabla^{M} f\right|^{2}+f^{2} H^{2}}
\end{align*}
on $M$ with Neumann boundary condition $\left\langle\nabla^{M} \psi
, \vec{n}\right\rangle=1$ on $\partial M$. Here, $\vec{n}$ denotes
the co-normal to $M$. Notice that according to the standard theory
of elliptic regularity, $\psi$ belongs to the class $C^{2, \beta}$
for each $0<\beta<1$;

The third fold is to derive the following inequality
\begin{align}\label{1.3}
    \int_{M} f^{\frac{n}{n-1}} d\mu \geq |B^{n}|
\end{align}
which is reduced from the original inequality (\ref{1.2}).

In addition to the Michael-Simon inequality for mean curvature, it
is certainly meaningful to establish the Michael-Simon type
inequality for $k$th mean curvature $\sigma_{k}(\kappa)$. In
\cite{CW2013}, Sun-Yung Alice Chang and Yi Wang presented the
following inequality in $\mathbb{R}^{n+1}$.

\begin{theorem}{(\cite{CW2013})}\label{th1.3}
    Let $i: M^{n} \rightarrow \mathbb{R}^{n+1}$ be an isometric immersion. Let $U$ be an open subset of $M$ and $u \in C_{c}^{\infty}(U)$ be a nonnegative function. For $m =2, \ldots, n-1$, if $M$ is $(m+1)$ -convex, then there exists a constant $C$ depending only on $n$ and $m$, such that for $1 \leqslant k \leqslant m$
    \begin{align*}
    &\left(\int_{M} \sigma_{k-1}(\kappa) u^{\frac{n-k+1}{n-k}} d \mu_{M}\right)^{\frac{n-k}{n-k+1}} \leq C \int_{M}\left(\sigma_{k}(\kappa) u+\sigma_{k-1}(\kappa)|\nabla u|+\cdots+\left|\nabla^{k} u\right|\right) d \mu_{M}
    \end{align*}
    If $m=n$, then the inequality holds when $M$ is $n$-convex. If $m=1$, then the inequality holds when $M$ is 1-convex. ($m=1$ case is a corollary of the Michael-Simon inequality).
\end{theorem}
In \cite{CW2013}, the method used to prove the above result is
optimal transport theory. M. Gromov first proposed the idea of
proving geometric inequalities by means of mappings between domains
and spheres which are optimal transport maps for special cases (see
e.g. \cite{M80,XZ17}). But what is the case in the equality holds?
And what is the sharp constant?

Motivated by Theorem \ref{th1.2} and Theorem \ref{th1.3}, we propose
the following conjecture about the sharp Michael-Simon inequality
for $k$th mean curvature in Euclidean space $\mathbb{R}^{n+1}$.
\begin{conjecture}\label{con1.4}
    Let $M$ be a compact hypersurface in $\mathbb{R}^{n+1}$ (possibly with boundary $\partial M$) and let $f$ be a positive smooth function on $M$. For $1\leq k\leq n-1$, there holds
    \begin{align}\label{1.4}
        \int_{M} \sqrt{\sigma_{k}^{2}f^{2}+\sigma_{k-1}^{2}|\nabla^{M} f|^{2}} +
        \int_{\partial M} \sigma_{k-1}f \geq n|B^{n}|^{\frac{1}{n-k+1}}\left(\int_{M}\sigma_{k-1}f^{\frac{n-k+1}{n-k}}\right)^{\frac{n-k}{n-k+1}}
    \end{align}
    where $\sigma_{k}:=\sigma_{k}(\kappa)$ is the $k$th mean curvature. Equality holds in (\ref{1.4}) if and only if $M$ is a sphere and $f$ is constant.
\end{conjecture}

We know that the constrained curvature flow is a very effective tool
for proving new and old geometric inequalities that do not generate
singularities during the flow. The first example of such flow was
proposed by Huisken \cite{Hui87} in Euclidean space
$\mathbb{R}^{n+1}$ and it is called volume preserving mean curvature
flow. Huisken proved that if the initial hypersurface is uniformly
convex, then the flow exists for all time and converges to a sphere
in $C^{\infty}$-topology as $t\to +\infty $. Notice that the
velocity function of this flow contains a global constrained term so
that the volume of the closed domain enclosed by the hypersurface is
unchanging and the area of the hypersurface is monotonically
decreasing. This property makes the evolving hypersurface converge
to the solution of the isoperimetric type problem in
$\mathbb{R}^{n+1}$. For similar curvature flows, one can refer to
\cite{And01,CS10,McC03,WangXia14} and the references therein for
details.

The globally constrained curvature flow is very difficult to compute
a priori estimates, mainly due to the global constraint term. In
\cite{GL2015}, Guan and Li defined a locally constrained curvature
flow in the space form, which is based on Minkowski identity. This
type of flow only requires the initial hypersurface to be a
starshaped to obtain the longtime existence and convergence.
Subsequently, locally constrained curvature flow attracted a lot of
attention, such as in \cite{GL09,HLW20,HuLi21,SXia19,WeiXiong21}.
However, its applications are mainly focused on proving the
isoperimetric type inequalities.

Mean curvature flow also has promising applications in physics. It
can not only describe the evolution of interfaces, such as
propagation at material interfaces, fluid free boundary motion, and
crystal growth (see e.g. \cite{WM1956,PS1997}), but is also widely
developed in fields such as image processing, computer-aided design,
and algorithms.

In this paper, we will introduce the new mean curvature type flow
and apply it to prove some remarkable inequalities such as
Michael-Simon inequality.

In the first part of the present paper, we introduce a new locally
constrained mean curvature flow $X:M \times [0,T) \to
\mathbb{R}^{n+1}$,
 $n\geq 2$ of starshaped hypersurfaces in $\mathbb{R}^{n+1}$ which satisfies
\begin{equation}\label{1.5}
    \begin{cases}
        \frac{\partial }{\partial t}X(x,t)=-\left(fH + \frac{n}{n-1}\frac{\partial f}{\partial X}\right)\nu(x,t)\\
        X(\cdot,0)=X_{0}(\cdot)
     \end{cases}
\end{equation}
where $\nu(x,t)$ and $H$ are the unit outer normal and the mean
curvature of $M_{t}=X(M,t)$, respectively, and $f: M \times[0,T,)
\to \mathbb{R}$ is a positive smooth function.

Suppose $M_{t}$ is starshaped with respect to a point $p$,  given by
a smooth embedding $X(\cdot, t): \mathbb{S}^{n} \to  M_{t} \subset
\mathbb{R}^{n+1}$. The radial function $r:\mathbb{S}^{n} \times
[0,T) \to \mathbb{R}_{+}$ represents the distance from $X(\cdot, t)$
to $p$, then we have $X(\xi,t)=r(\xi,t)\xi, \xi \in \mathbb{S}^{n}$
and $f$ can be expressed as $f\left[r\left(\xi,t\right)\right]$.
Therefore, the evolution problem (\ref{1.5}) can be expressed, in
terms of $r$, as a scalar PDE
\begin{equation}\label{1.6}
    \begin{cases}
        \frac{\partial}{\partial t}r=-\left(fH+\frac{n}{n-1} \frac{\partial f}{\partial r}\sqrt{1+r^{-2}|\nabla r|^{2}}\right)\sqrt{1+r^{-2}|\nabla r|^{2}}
        \\
        r(\cdot,0)=r_{0}
    \end{cases}
\end{equation}
we will prove that the flow (\ref{1.5}) has longtime existence and
convergence provided that $f$ satisfies the following assumption.
\begin{assumption}\label{as1.5}
    $\widehat{f}(r):=(n-1)\frac{1}{r^{2}}f(r)+\frac{1}{r}\frac{\partial f}{\partial r}(r)$ is monotonically increasing with respect to $r$, where $r>0$, and there exists a zero point for $\widehat{f}(r)$.
\end{assumption}
Now we state the main properties of flow (\ref{1.5}).
\begin{theorem}\label{th1.6}
    Let $X_{0}:M \to \mathbb{R}^{n+1}(n\geq 2)$ be a smooth embedding of a compact hypersurface $M$ in $\mathbb{R}^{n+1}$ such that $M_{0}=X_{0}(M)$ is starshaped, and assume that $f$ satisfies Assumption \ref{as1.5}. Then the flow (\ref{1.5}) has a unique smooth solution $M_{t}=X_{t}(M)$ for all time $t\in [0,+\infty)$. Moreover, $M_{t}$ converges to a sphere as $t\to +\infty $ in $C^{\infty}$-topology.
\end{theorem}
In addition, we can apply the convergence result of flow (\ref{1.5})
to prove the inequality (\ref{1.2}) for starshaped hypersurfaces.
\begin{theorem}\label{th1.7}
    Let $M$ be a smooth, compact, starshaped hypersurface in $\mathbb{R}^{n+1}$ (possibly with boundary $\partial M$ ). For a positive function $f \in C^{\infty}(M)$, there holds
     \begin{align*}
         \int_M \sqrt{|\nabla^{M} f|^{2}+f^{2}H^{2}}d\mu + \int_{\partial M}f d\mu \geq n|B^{n}|^{\frac{1}{n}}\left(\int_{M} f^{\frac{n}{n-1}} d\mu\right)^{\frac{n-1}{n}}
     \end{align*}
    Equality holds if and only if $M$ is a sphere and $f$ is constant.
\end{theorem}
Further, we briefly explain the proof of Theorem \ref{th1.7} in
three steps.

Step 1. We apply the divergence theorem and the standard elliptic
regular theory to simplify the original inequality (\ref{1.2}) by
scaling;

Step 2. We will show that $\int_{M_{t}}f^{\frac{n}{n-1}}d\mu_{t}$
decreases monotonically along the flow (\ref{1.5}), this property
makes the inequality (\ref{1.2}) hold for any evolving
hypersurfaces;

Step 3. Following the convergence result of flow (\ref{1.5}), we
obtain the sharp constant in the inequality (\ref{1.2}) and the
necessary and sufficient condition for the establishment of the
equality.

In the second part of this paper, we will introduce a new mean
curvature type flow and use it to solve  Conjecture \ref{con1.4}.
Let $X_{0}:M \to \mathbb{R}^{n+1}$ be a smooth embedding such that
$M$ is a closed, convex hypersurface in Euclidean space
$\mathbb{R}^{n+1}$. We consider a smooth family of embeddings $X:M
\times [0,T) \to \mathbb{R}^{n+1}$ satisfying
\begin{equation}\label{1.7}
    \begin{cases}
        \frac{\partial}{\partial t}X(x,t)=\left(1-h\frac{E_{k}(\kappa)}{E_{k-1}(\kappa)}\right)\nu(x,t),\quad k=1,\cdots ,n\\
        X(\cdot,0)=X_{0}(\cdot)
    \end{cases}
\end{equation}
where $\nu(x,t)$ and $\kappa=(\kappa_{1},\cdots,\kappa_{n})$ are the
unit outer normal and the principal curvatures of $M_{t}=X(M,t)$
respectively, $E_{k}(\kappa)=\binom{n}{k}^{-1}\sigma_{k}(\kappa)$ is
the normalized $k$th mean curvature and $h(\nu)=\langle
X,\nu\rangle$ is the support function of $M_{t}$. When $k=1$, by
scaling, the flow (\ref{1.7}) is equivalent to the mean curvature
type flow $\frac{\partial }{\partial t}X=\left(n-Hh \right)\nu$
introduced in \cite{GL2015}.

Our construction of the flow (\ref{1.7}) is inspired by the fully
nonlinear flow defined by Guan and Li in \cite{GL2018}, which takes
the following form
\begin{align}\label{1.8}
    \frac{\partial }{\partial t}X(x,t)&=\left(c_{k-1}-h\frac{\sigma_{k}(\kappa)}{\sigma_{k-1}(\kappa)}\right)\nu(x,t)=c_{k-1}\left(1-h\frac{E_{k}(\kappa)}{E_{k-1}(\kappa)}\right)\nu(x,t)
\end{align}
where
$c_{k-1}=\frac{\sigma_{k}(I)}{\sigma_{k-1}(I)}=\frac{n+1-k}{k}$.
They show that starting from a smooth, closed, convex hypersurface
in $\mathbb{R}^{n+1}(n \geq 2)$, the solution of the flow
(\ref{1.8}) exists for all positive time and converges smoothly and
exponentially to a sphere.  A nice characteristic of  flow
(\ref{1.8}) is that the $(k-1)$th quermassintergral is unchanging
and the $k$th quermassintergral is monotonically decreasing, so it
can be used to prove the Alexandrov-Fenchel inequality.
\begin{definition}\label{def1.8}
    For a bounded domain $\Omega$ in $\mathbb{R}^{n+1}$ with smooth boundary $M=\partial\Omega$, it is called static convex if its second fundamental form satisfies
    \begin{align}\label{1.9}
        h_{ij}\geq h^{-1}g_{ij}>0
    \end{align}
    everywhere on $M$, where $h^{-1}$ is the inverse of the support function.
\end{definition}
Concerning the flow (\ref{1.7}), we obtain a new property that if
the initial hypersurface is static convex, then the solution of the
flow (\ref{1.7}) remains static convex for all $t>0$. It is worth
noting that this property is crucial in proving Conjecture
\ref{1.4}. From the inequality (\ref{1.9}), we get the static
convexity implies the strict convexity. Further, it can be found
that the velocity functions of flow (\ref{1.7})  and flow
(\ref{1.8}) differ only by a constant multiple. Therefore, the
longtime existence and convergence of flow (\ref{1.7}) can be
obtained when the initial hypersurface is static convex.
\begin{theorem}\label{th1.9}
   Let $X_{0}:M \to \mathbb{R}^{n+1}(n\geq 2)$ be a smooth embedding of a closed, static convex hypersurface $M_{0}=X_{0}(M)$ in $\mathbb{R}^{n+1}$. Then the flow (\ref{1.7}) has a unique smooth solution $M_{t}=X_{t}(M))$ for all time $t\in[0,+\infty)$. Moreover, $M_{t}$ is static convex for each $t\geq 0$ and it converges exponentialy to a sphere of radius $R$ in $C^{\infty}$-topology as $t \to +\infty$, where the radius $R$ determined by $V_{k-1}(\Omega_{0})=V_{k-1}(B^{n+1}_{R})$.
\end{theorem}
Inspired by the idea of the proof of the Michael-Simon inequality
for mean curvature, We can use the properties of flow (\ref{1.7}) to
prove a new sharp Michael-Simon inequality for $k$th mean curvature
in $\mathbb{R}^{n+1}$. To facilitate our discussion, let $M_{t}$ be
parametrized by the inverse Gauss map $X_{t}:\mathbb{S}^{n} \to
M_{t} \subset  \mathbb{R}^{n+1}$, then the positive function $f \in
C^{\infty}(M_{t})$ can be expressed as
$f\left[h\left(\nu\right)\right]$. Before presenting the detailed
results, we need to make the following assumption about $f$.
\begin{assumption}\label{as1.10}
    $g(h):=f^{\frac{n-k+1}{n-k}}(h)$ is convex and monotonically increasing with respect to $h$, where $h>0$.
\end{assumption}

\begin{theorem}\label{th1.11}
Let $M$ be a smooth, compact and static convex hypersurface in
$\mathbb{R}^{n+1}$ (possibly with boundary $\partial M$), and
$\Omega$ be the domain enclosed by $M$. Assume that $f$ satisfies
Assumption \ref{as1.10}. Then for any $1\leq k \leq n-1$, there
holds
    \begin{align}\label{1.10}
        \int_{M} & \sqrt{\sigma_{k}^{2}f^{2}+\sigma_{k-1}^{2}|\nabla^{M} f|^{2}}d\mu +\int_{\partial M} \sigma_{k-1}f d\mu \notag \\ 
        &\geq n\left(y_{k}\circ z_{k-1}^{-1}(V_{k-1}(\Omega))\right)^{\frac{1}{n+1-k}}\left(\int_{M}\sigma_{k-1}f^{\frac{n+1-k}{n-k}}d\mu\right)^{\frac{n-k}{n+1-k}}
    \end{align}
    where $y_{k}(r)=\binom{n}{k}f^{\frac{n+1-k}{n-k}}z_{k}(r)$, $z_{k}(r)=V_{k}(B^{n+1}_{r})$, the $k$th quermassintergral for the sphere of radius $r$,
    and $z_{k}^{-1}$ is the inverse function of $z_{k}$. Equality holds in (\ref{1.10}) if and only if $M$ is a sphere and $f$ is constant.
\end{theorem}

In particular, if $M$ is a sphere with radius $R$ determined by
$V_{k-1}(\Omega)=V_{k-1}(B^{n+1}_{R})$ and we consider the special
case that $f=R^{-(n-k)}$ in Theorem \ref{th1.11} and derive the
following conclusion.
\begin{corollary}\label{col1.12}
    Let $M$ and $f$ be given as above, then the inequality (\ref{1.4}) for static convex hypersurface is established. Equality holds in (\ref{1.4}) if and only if $M$ is a sphere and $f$ is constant.
\end{corollary}

\begin{remark}\label{rem1.13}
    (1) Corollary \ref{col1.12} is the weaker form of Conjecture \ref{con1.4}.

    (2) If $k=1$, the inequality (\ref{1.4}) is the sharp Michael-Simon inequality for mean curvature.
\end{remark}
The paper is organized as follows. In Section 2, we recall some
geometries of hypersurfaces in Euclidean space and collect some
basic properties of normalized elementary symmetric functions. In
Section 3, we derive the evolution equations along the flows
(\ref{2.6}) and (\ref{1.7}). In Section 4, by calculating the
required prior estimates, we get the longtime existence of the flow
(\ref{1.5}), and the convergence of flow (\ref{1.5}) is obtained by
a detailed estimation of the gradient of the radial function. The
new proof of the sharp Michael-Simon inequality (\ref{1.2}) will be
given in Section 5. In Section 6, we apply the tensor maximum
principle to show that the flow (\ref{1.7}) preserves the static
convexity of the evolving hypersurfaces, and then complete the proof
of Theorem \ref{th1.9}. The convergence result of the flow
(\ref{1.7}) will be used in Section 7 to prove inequality
(\ref{1.4}) for static convex hypersurface.

\section{\bf Preliminaries}

In this section, we first collect three parametrizations of
(locally) hypersurfaces embedded in $\mathbb{R}^{n+1}$. These are
parametrization by radial function, by a graph of a function, by
support function. The first two parameterizations are used for the
flow (\ref{1.5}), where the second parameterization is used in
estimating the principal curvatures. the last one is used for the
flow (\ref{1.7}). Last, we review the properties of normalized
elementary symmetric functions.

\noindent{\bf 2.1~~ Parametrization by radial graph}

\noindent As $M\subset \mathbb{R}^{n+1}$ is a smooth, compact,
starshaped hypersurface with respect to a point, we may suppose to
be the origin of $\mathbb{R}^{n+1}$. $M$ can be represented as
\begin{align*}
    M=\{r(\xi)\xi: \xi \in \mathbb{S}^{n}\}
\end{align*}
where $r=|X|$ is the radial function and $\xi=\frac{X}{|X|}$. Let
$\left\{e_{1},\cdots,e_{n}\right\}$ be a local orthonormal
coordinate system of $\mathbb{S}^{n}$, and denote  $\nabla$ be the
gradient on $\mathbb{S}^{n}$. The induced metric on $M$ now has the
form (see e.g.\cite{G})
\begin{align*}
    g_{ij}=\left\langle \nabla_{i}X, \nabla_{j}X\right\rangle=r^{2}e_{ij}+\nabla_{i}r\nabla_{j}r
\end{align*}
where $\left\langle \cdot,\cdot \right\rangle $ denotes the standard
inner product on $\mathbb{R}^{n+1}$ and $e_{ij} $ is the metric of
$\mathbb{S}^{n}$. The inverse metric $\left(g_{ij}\right)^{-1}$ is
\begin{align*}
    g^{ij}=r^{-2}(e^{ij}-\frac{\nabla^{i}r\nabla^{j}r}{r^{2}+{|{\nabla r}|}^{2}})
\end{align*}
where $\nabla^{i}=e^{ik}\nabla_{k}$. The outer unit normal is given
by
\begin{align}\label{2.1}
    \nu=(r \xi -\nabla r)\left(r^{2}+|\nabla r|^{2}\right)^{-\frac{1}{2}}
\end{align}
The second fundamental form and the mean curvature of $M$ are as
follows
\begin{align}\label{2.2}
    h_{i j}&=-\left\langle\nabla_{i j} X, \nu\right\rangle=\left(r^{2}+|\nabla r|^{2}\right)^{\frac{1}{2}}\left(r^{2} e_{i j}+2 \nabla_{i} r \nabla_{j}r - r\nabla_{i j}r\right)
\end{align}
\begin{align}\label{2.3}
    H&=\sum g^{ij}h_{ij}=\frac{1}{r\sqrt{1+r^{-2}|\nabla r|^{2}}}\left[n-\left(e^{ij}-\frac{\nabla^{i}\omega \nabla^{j}\omega}{1+|\nabla \omega|^{2}}\right)\nabla_{ij}\omega\right]
\end{align}
where $\nabla_{i j}=\nabla_{i} \nabla_{j}$ and $\omega=lnr$. The
principal curvatures $\kappa=\left( \kappa_{1},\cdots, \kappa_{n}
\right)$ of $M$ are the eigenvalues of the second fundamental form
relative to the metric, or equivalently, they are the solutions of
the equation
\begin{align*}
    \operatorname{det}\left(\lambda_{ij}-\kappa \delta_{ij}\right)=0
\end{align*}
where the symmetric matrix $\left[\lambda_{ij}\right]$ is given by
(see e.g.\cite{Urbas1990})
\begin{align}\label{2.4}
    \left[\lambda_{ij}\right]=\left[g^{ik}\right]^{\frac{1}{2}}\left[h_{kl}\right]\left[g^{lj}\right]^{\frac{1}{2}}
\end{align}
$\left[g^{i j}\right]^{\frac{1}{2}}$ denotes the positive square
root of $\left[g^{i j}\right]$ and can be expressed as
\begin{align}\label{2.5}
    \left[g^{ij}\right]^{\frac{1}{2}}=r^{-1}\left[e^{ij}-\frac{\nabla^{i}r \nabla^{j}r}{\sqrt{r^{2}+|\nabla r|^{2}}\left(r^{2}+\sqrt{r^{2}+|\nabla r|^{2}}\right)}\right]
\end{align}

Note that, according to \cite{Urbas1990}, if the flow equation is
\begin{align}\label{2.6}
   \frac{\partial}{\partial t}X(x,t)=\Phi(x,t) \nu(x,t)
\end{align}
where $\Phi(x,t)$ is a smooth function on $M_{t}$, then the radial
function $r$ satisfies the following equation
\begin{align}\label{2.7}
    \begin{cases}
        \frac{\partial}{\partial t}r = \sqrt{1+r^{-2}|\nabla r|^{2}}\Phi(x,t) \qquad \text{on}\quad\mathbb{S}^{n} \times \mathbb{R}_{+}
       \\
       r(\cdot,0)=r_{0}
    \end{cases}
\end{align}
Thus the flow (\ref{1.5}) can be converted into (\ref{1.6}).
\\
\noindent{\bf 2.2~~ Parametrization by local graph }

\noindent It is known that $M$ is smooth, and by the implicit
function theorem, we may assume that $M$ can be locally represented
by a graph of a function on hyperplane $\Sigma \subset
\mathbb{R}^{n+1}$. By rotating coordinates, $\Sigma$ can become
$\mathbb{R}^{n}$, then there is a $C^{\infty}$ function $u:\Sigma
\to \mathbb{R} $ defined in a neighborhood $U$ of a point $P \in M$
that satisfies
\begin{align*}
    M \cap  U  =\{(\eta,u(\eta)); \eta \in U \subset  \Sigma\}
\end{align*}
The metric, the inverse metric, the second fundamental form and the
principal curvatures of $M$ are given by (see e.g.\cite{Urbas1990})
\begin{align*}
 g_{i j}=\delta_{i j}+D_{i} u D_{j} u ;\quad
 g^{i j}=\delta_{i j}-\frac{D_{i} u D_{j} u}{1+|D u|^{2}};\quad
 h_{i j}=\frac{D_{i j} u}{\left(1+|D u|^{2}\right)^{\frac{1}{2}}}
\end{align*}
and
\begin{align}\label{2.8}
    \lambda_{ij}=\frac{1}{v}\left\{D_{ij}u -\frac{D_{i}uD_{l} u D_{jl} u}{v(1+v)}-\frac{D_{j}u D_{l}uD_{il} u}{v(1+v)} +\frac{D_{i} u D_{j} u D_{k} u D_{l} u D_{k l} u}{v^{2}(1+v)^{2}}\right\}
\end{align}
where $D$ is the usual gradient in $\mathbb{R}^{n}$, and
$v=\left(1+|D u|^{2}\right)^{\frac{1}{2}}$ (see \cite{CNS}).
Together with (\ref{2.6}), we can obtain that $u(\cdot,t):\Sigma
\times [0,T) \to \mathbb{R}$ is the solution of the following
equation
\begin{align}\label{2.9}
    \begin{cases}
       D_{t} u=- \sqrt{ \left(1+|D u|^{2}\right) } \Phi(x,t) \qquad \text{on}\quad\Sigma \times \mathbb{R}_{+}
       \\
       u(\cdot,0)=u_{0}
    \end{cases}
\end{align}
\\
\noindent{\bf 2.3~~Parametrization by support function}

\noindent Let $M$ be a smooth, compact, strictly convex hypersurface
in $\mathbb{R}^{n+1}$, one may suppose that the origin in its
interior, then $M$ can be parametrized by the inverse Gauss map $X:
\mathbb{S}^{n} \rightarrow M \subset \mathbb{R}^{n+1}$, defined as
\begin{align*}
    X(\nu)=h(\nu) \nu+\nabla h(\nu)
 \end{align*}
where $\nu$ is the unit outer normal and $h(\nu)$ is the support
function. The second fundamental form of $M$ is given by (see e.g.
\cite{Urbas1991})
\begin{align}\label{2.10}
   h_{i j}=\nabla_{j}\nabla_{i} h+h e_{i j}
\end{align}
From the Guass-Weingarten formula
$\nabla_{i}\nu=h_{ik}g^{kl}\nabla_{l}X$, we have
\begin{align}\label{2.11}
   e_{i j}=\left\langle\nabla_{i}\nu,\nabla_{j}\nu\right\rangle=h_{i k} g^{k l} h_{j l}
\end{align}
and
\begin{align*}
    g_{ij}=h_{ik}e^{kl}h_{jl}
\end{align*}
According to (\ref{2.6}), the general evolution equation of $h$ is
\begin{align*}
    \begin{cases}
        \frac{\partial}{\partial t}h = \Phi(x,t) \qquad \text{on} \quad \mathbb{S}^{n} \times \mathbb{R}_{+}\\
        h(\cdot,0)=h_{0}
    \end{cases}
\end{align*}
Thus, the flow equation (\ref{1.7}) can be expressed in terms of $h$
as the following equation
\begin{align}\label{2.12}
    \begin{cases}
        \frac{\partial}{\partial t}h=1-h\frac{E_{k}(\kappa)}{E_{k-1}(\kappa)} \qquad \text{on} \quad \mathbb{S}^{n} \times \mathbb{R}_{+}\\
        h(\cdot,0)=h_{0}
    \end{cases}
\end{align}

\begin{lemma}
   Let $M$ be a smooth, strictly convex hypersurface in $\mathbb{R}^{n+1}$, and the position vecter $X:\mathbb{S}^{n} \to \mathbb{R}^{n+1}$. Then the support function $h$ satisfies
   \begin{align}
    &(1)\nabla_{i}h=<X,x_{k}>h^{k}_{i}\label{2.13}\\
    &(2)\nabla_{j}\nabla_{i}h=h_{ij}-h(h^{2})_{ij} \label{2.14}
   \end{align}
   where $\{x_{1},\cdots ,x_{n}\}$ is a coordinate system in the tangent space of $M$.
\end{lemma}
\noindent{\it \bf Proof.}~~(1) Differentiating $h=\langle X,
\nu\rangle$ gives
\begin{align*}
    \nabla_{i}h=\nabla_{i}\langle X, \nu\rangle=\langle \nabla_{i}X, \nu\rangle +\langle X, \nabla_{i}\nu\rangle
    =\langle X, \nabla_{i}\nu\rangle
\end{align*}
and combining the Guass-Weingarten formula, we have
\begin{align*}
    \nabla_{i}h=h^{i}_{l}\langle X, \nabla_{l}X\rangle =h^{i}_{l}\langle X, x_{l}\rangle
\end{align*}

(2) The third equation can be directly obtained by (\ref{2.10}).
\hfill${\square}$

\noindent{\bf 2.4~~ Normalized elementary symmetric functions }

\noindent For each $k=1, \ldots, n$, the normalized $k$th elementary
symmetric functions for $\kappa=\left(\kappa_{1}, \ldots,
\kappa_{n}\right)$ are
\begin{align*}
    E_{k}(\kappa)=\binom{n}{k}^{-1} \sigma_{k}(\kappa)=\binom{n}{k}^{-1} \sum_{1 \leq i_{1}<\ldots<i_{k} \leq n} \kappa_{i_{1}} \cdots \kappa_{i_{k}}
\end{align*}
and we can set $E_{0}(\kappa)=1$ and $E_{k}(\kappa)=0$ for $k>n$. If
$A=[A_{ij}] $ is an $n \times n$ symmetric matrix and
$\kappa=\kappa(A)=(\kappa_{1},\cdots,\kappa_{n})$ are the
eigenvalues of $A$, then $E_{k}(A)=E_{k}(\kappa(A))$ can be
expressed as
\begin{align*}
    E_{k}(A)=\frac{(n-k)!}{n!} \delta_{i_{1} \ldots i_{k}}^{j_{1} \ldots j_{k}} A_{i_{1} j_{1}} \cdots A_{i_{k} j_{k}}, \quad k=1, \ldots, n
\end{align*}
Now, let's review some of the properties of the normalized $k$th
elementary symmetric functions (see e.g.\cite{G2013}).
\begin{lemma}
   Let $\dot{E}_{k}^{i j}=\frac{\partial E_{k}}{\partial A_{i j}}$, then we have
   \begin{align}
    \sum_{i, j} \dot{E}_{k}^{i j} g_{ij} &=k E_{k-1} \label{2.15}\\
    \sum_{i, j} \dot{E}_{k}^{i j} A_{i j} &=k E_{k} \label{2.16}\\
    \sum_{i, j} \dot{E}_{k}^{i j}\left(A^{2}\right)_{i j} &=n E_{1} E_{k}-(n-k) E_{k+1} \label{2.17}
    \end{align}
   where $\left(A^{2}\right)_{i j}=\sum_{k=1}^{n} A_{i k} A_{k j}$
\end{lemma}
Next, we recall the well-know Minkowski identity (see
e.g.\cite{GL2015}).
\begin{lemma}   Let $M$ be a smooth closed hypersurface in $\mathbb{R}^{n+1}$. Then
   \begin{align}\label{2.18}
   \int_{M} E_{k-1}(\kappa) d\mu= \int_{M} hE_{k}(\kappa) d\mu
   \end{align}
\end{lemma}
\begin{lemma}
    If $\kappa \in \Gamma_{k}^{+}$, the following inequality is called the Newton-MacLaurin inequality
   \begin{align}\label{2.19}
      E_{m+1}(\kappa) E_{k-1}(\kappa) \leq E_{k}(\kappa) E_{m}(\kappa), \quad 1 \leq k \leq m
   \end{align}
   where $\Gamma ^{+}_{k}=\left\{x\in \mathbb{R}^{n}:E_{i}(x)>0, i=1,\cdots,k\right\}$. Equality holds if and only if $\kappa_{1}=\cdots=\kappa_{n}$.
\end{lemma}
Let $F(\kappa)=\frac{E_{k}(\kappa)}{E_{k-1}(\kappa)}$, $F(\kappa)$
is a smooth symmetric function on $\mathbb{R}^{n}$, $\dot{F}^{p}$
and $\ddot{F}^{pq}$ are the components of the first and second
derivatives of $F$ with respect to its argument, i.e.
\begin{align*}
    \dot{F}^{p}(\kappa)=\frac{\partial F(\kappa)}{\partial \kappa_{p}};\quad \ddot{F}^{pq}(\kappa)=\frac{\partial^{2} F(\kappa)}{\partial \kappa_{p} \partial \kappa_{q}}
\end{align*}
For a diagonal matrix $A$ with eigenvalues $\kappa=\kappa(A)$,
similarly, we can view $F(\kappa)$ as a smooth symmetric function
$F(A)=F(\kappa(A))$, where  $F(A)=\frac{E_{k}(A)}{E_{k-1}(A)}$. In a
local orthonornal frame, the first and second derivatives of $F(A)$
satisfy (see \cite{And94},  \cite{And07})
\begin{align*}
    \dot{F}^{pq}(A)=\dot{F}^{p}(\kappa)\delta_{pq}
\end{align*}
and
\begin{align*}
    \ddot{F}^{pq,ij}(A)B_{pq}B_{ij}=\sum_{p, q}\ddot{F}(\kappa)B_{pp}B_{qq}+2\sum_{p<q}\frac{\dot{F}^{p}(\kappa)-\dot{F}^{q}(\kappa)}{\kappa_{p}-\kappa_{q}}(B_{pq})^{2}
\end{align*}
where $B\in Sym(n)$. The later formula makes sense as a limit in the
case of $\kappa_{p}=\kappa_{q}$. Using (\ref{2.15}), (\ref{2.17})
and (\ref{2.19}), we have the following corollary.
\begin{corollary}
   Let $F(A)=\frac{E_{k}(A)}{E_{k-1}(A)}$ and $\kappa(A) \in \Gamma_{k}^{+}$, then
   \begin{align}
      &1 \leq \sum_{i, j} \dot{F}^{i j} g_{ij} \leq k \label{2.20}\\
      &F^{2} \leq \sum_{i, j} \dot{F}^{i j}\left(A^{2}\right)_{i j} \leq(n-k+1) F^{2} \label{2.21}
   \end{align}
\end{corollary}
\begin{lemma}\label{le2.6}
    See\cite{And07}, $F(\kappa)=\frac{E_{k}(\kappa)}{E_{k-1}(\kappa)}$ satisfies the following, where $\kappa \in \Gamma_{k}^{+}$\\
    (1) $F(\kappa)$ is strictly increasing, i.e. $\dot{F}^{p}=\frac{\partial F}{\partial \kappa_{q}}>0$ on $\Gamma_{k}^{+}$, $\forall p=1,\cdots,n$;\\
    (2) $F(\kappa)$ is homogeneous of degree 1, i.e. $F(a\kappa)=aF(\kappa)$ for any $a>0$;\\
    (3) $F(\kappa)$ is strictly positive on $\Gamma_{k}^{+}$ and is normalized such that $F(I)=1$;\\
    (4) $F(\kappa)$ is concave;\\
    (5) $F(\kappa)$ is inverse concave, i.e. the function
    \begin{align*}
        F_{*}(\kappa_{1},\cdots,\kappa_{n})=F(\kappa_{1}^{-1},\cdots,\kappa_{n}^{-1})^{-1}
    \end{align*}
    is concave.
\end{lemma}
\section{\bf Evolution equation}
\noindent Along the general flow (\ref{2.6}) that
\begin{align*}
     \frac{\partial}{\partial t}X(x,t)=\Phi(x,t) \nu(x,t)
\end{align*}
in Euclidean space $\mathbb{R}^{n+1}$, we have the following evolution equations (see \cite{GL09}).
\begin{lemma}
    \begin{align}
        \frac{\partial}{\partial t}g_{ij}&=2\Phi h_{ij}\label{3.1}\\
        \frac{\partial}{\partial t}d\mu_{t}&=nE_{1}\Phi d\mu_{t}\label{3.2}\\
        \frac{\partial}{\partial t}h_{ij}&=-\nabla_{j}\nabla_{i}\Phi+\Phi(h^{2})_{ij}\label{3.3}\\
        \frac{\partial}{\partial t}h^{j}_{i}&=-\nabla^{j}\nabla_{i}\Phi-\Phi(h^{2})^{j}_{i}\label{3.4}\\
        \frac{\partial}{\partial t}E_{k-1}&=\frac{\partial E_{k-1}}{\partial h^{j}_{i}}\frac{\partial h^{j}_{i}}{\partial t}=\dot{E}^{ij}_{k-1}\left(-\nabla_{j}\nabla_{i}\Phi - \Phi(h^{2})_{ij} \right) \label{3.5}
    \end{align}
\end{lemma}

\begin{lemma}
    Along the flow (\ref{1.7}), we have the following evolution equations.\\
    (1) The second fundamental form evolves
    \begin{align}
        \frac{\partial}{\partial t}h_{ij}=&h\dot{F}^{kl}\nabla_{k}\nabla_{l}h_{ij}+h\ddot{F}^{kl,pq}\nabla_{i}h_{kl}\nabla_{j}h_{pq}+\nabla_{i}F\nabla_{j}h+\nabla_{i}h\nabla_{j}F \notag\\
        &+\left(F+h\dot{F}^{kl}(h^{2})_{kl}\right)h_{ij}+(1-3hF)(h^{2})_{ij}\label{3.6}
    \end{align}
    (2) Let $S_{ij}=h_{ij}-h^{-1}g_{ij}$, it evolves as
    \begin{align}
        \frac{\partial}{\partial t}S_{ij}=&h\dot{F}^{kl}\nabla_{k}\nabla_{l}S_{ij}+h^{-2}\dot{F}^{kl}\nabla_{k}h\nabla_{l}hg_{ij}+h\ddot{F}^{kl,pq}\nabla_{i}h_{kl}\nabla_{j}h_{pq}+\nabla_{i}h\nabla_{j}F+\nabla_{j}h\nabla_{j}F  \notag \\
        &+(1-3hF)(S^{2})_{ij}+\left(h\dot{F}^{kl}(h^{2})_{kl}-3F\right)S_{ij}      +2\dot{F}^{kl}(h^{2})_{kl}g_{ij}-2h^{-1}Fg_{ij} \label{3.7}
    \end{align}
\end{lemma}
\noindent{\it \bf Proof.}~~(1) We have know that
$\frac{\partial}{\partial
t}h_{ij}=-\nabla_{i}\nabla_{j}\Phi+\Phi(h^{2})_{ij}$ and
$\Phi=(1-hF)$, then
\begin{align}
    \frac{\partial}{\partial t}h_{ij}&=-\nabla_{i}\nabla_{j}(1-hF)+(1-hF)(h^{2})_{ij} \notag \\
    &=F\nabla_{i}\nabla_{j}h+h\nabla_{i}\nabla_{j}F+\nabla_{i}F\nabla_{j}h+\nabla_{i}h\nabla_{j}F+(1-hF)(h^{2})_{ij} \label{3.8}
\end{align}
where $(h^{2})_{ij}=h^{k}_{i}h_{kj}$. By the Simons's identity
$\nabla_{k}\nabla_{l}h_{ij}=\nabla_{i}\nabla_{j}h_{kl}-(h^{2})_{kl}h_{ij}+(h^{2})_{ij}h_{kl}$,
we have
\begin{align}
    \nabla_{i}\nabla_{j}F&=\nabla_{i}\left(\dot{F}^{kl}\nabla_{j}h_{kl}\right)=\ddot{F}^{kl,pq}\nabla_{i}h_{pq}\nabla_{j}h_{kl}+\dot{F}^{kl}\nabla_{i}\nabla_{j}h_{kl} \notag \\
    &=\dot{F}^{kl}\nabla_{k}\nabla_{l}h_{ij}+\ddot{F}^{kl,pq}\nabla_{i}h_{kl}\nabla_{j}h_{pq}+\dot{F}^{kl}(h^{2})_{kl}h_{ij}-\dot{F}^{kl}h_{kl}(h^{2})_{ij} \label{3.9}
\end{align}
Since $F$ is homogeneous of 1, then $F=\dot{F}^{kl}h_{kl}$.
Substituting (\ref{2.14}) and (\ref{3.9}) into (\ref{3.8}), we get
\begin{align*}
    \frac{\partial}{\partial t}h_{ij}=&h\dot{F}^{kl}\nabla_{k}\nabla_{l}h_{ij}+h\ddot{F}^{kl,pq}\nabla_{i}h_{kl}\nabla_{j}h_{pq}+\nabla_{i}F\nabla_{j}h+\nabla_{i}h\nabla_{j}F\\
    &+h\dot{F}^{kl}(h^{2})_{kl}h_{ij}+Fh_{ij}-hF(h^{2})_{ij}+(1-2hF)(h^{2})_{ij}\\
    =&h\dot{F}^{kl}\nabla_{k}\nabla_{l}h_{ij}+h\ddot{F}^{kl,pq}\nabla_{i}h_{kl}\nabla_{j}h_{pq}+\nabla_{i}F\nabla_{j}h+\nabla_{i}h\nabla_{j}F\\
    &+\left(F+h\dot{F}^{kl}(h^{2})_{kl}\right)h_{ij}+(1-3hF)(h^{2})_{ij}
\end{align*}
\\
(2) Differentiating $S_{ij}=h_{ij}-h^{-1}g_{ij}$ with respect to $t$
and substituting (\ref{2.12}), (\ref{3.1}) and (\ref{3.6}) yields
\begin{align*}
    \frac{\partial}{\partial t}S_{ij}=&h\dot{F}^{kl}\nabla_{k}\nabla_{l}h_{ij}+h\ddot{F}^{kl,pq}\nabla_{i}h_{kl}\nabla_{j}h_{pq}+\nabla_{i}F\nabla_{j}h+\nabla_{i}h\nabla_{j}F\\
    &+\left(F+h\dot{F}^{kl}(h^{2})_{kl}\right)h_{ij}+(1-3hF)(h^{2})_{ij}\\
    &+h^{-2}(1-hF)g_{ij}-2h^{-1}(1-hF)h_{ij}
\end{align*}
and by direct calculation, we get the following formula
\begin{align}
    \nabla_{k}\nabla_{l}S_{ij}&=\nabla_{k}\nabla_{l}h_{ij}+g_{ij}h^{-2}\nabla_{k}\nabla_{l}h-2h^{-3}g_{ij}\nabla_{k}h\nabla_{l}h\label{3.10}\\
    (S^{2})_{ij}&=S^{k}_{i}S_{kj}=(h^{2})_{ij}-2h^{-1}h_{ij}+h^{-2}g_{ij} \label{3.11}
\end{align}
Thus, we have
\begin{align*}
    \frac{\partial}{\partial t}S_{ij}= &h\dot{F}^{kl}\nabla_{k}\nabla_{l}S_{ij}+2h^{-2}g_{ij}\dot{F}^{kl}\nabla_{k}h\nabla_{l}h+h\ddot{F}^{kl,pq}\nabla_{i}h_{kl}\nabla_{j}h_{pq}+\nabla_{i}h\nabla_{j}F+\nabla_{j}h\nabla_{j}F\\
    &+(1-3hF)(S^{2})_{ij}+\left(h\dot{F}^{kl}(h^{2})_{kl}-3F\right)S_{ij}+2\dot{F}^{Kl}(h^{2})_{kl}g_{ij}-2h^{-1}Fg_{ij}
\end{align*}
\hfill${\square}$

\section{\bf Long time existence and smooth convergence of flow \ref{1.5}}
\noindent To obtain the longtime existence of the flow (\ref{1.5}),
we need to derive a priori estimates, thus using the standard theory
of parabolic partial differential equations. We already know that
the flow (\ref{1.5}) can be converted into scaler PDE (\ref{1.6}),
or equivalently,
\begin{align}\label{4.1}
    \begin{cases}
        \frac{\partial}{\partial t}r=-\sqrt{1+r^{-2}|\nabla r|^{2}}fH-\frac{n}{n-1} \frac{\partial f}{\partial r}(1+r^{-2}|\nabla r|^{2})
        \\
        r(\cdot,0)=r_{0}
    \end{cases}
\end{align}
Firstly, we perform $C^{0}$ estimate of $r$.
\begin{lemma}
    Let $r \in C^{\infty}\left(\mathbb{S}^{n} \times [0,T) \right)$ be a smooth solution to the initial value problem (\ref{4.1}), then there are positive constants $C_{1}=C_{1}( \delta, r_{min}(0))$ and $C_{2}=C_{2}( \varepsilon, r_{max}(0))$, such that
    \begin{align}\label{4.2}
          C_{1} \leq r(\cdot,t) \leq C_{2}
    \end{align}
\end{lemma}
\noindent{\it \bf Proof.}~~ Set $r_{min}(t):=\underset{x\in
\mathbb{S}^{n}}{min}  r(\cdot,t)$, we have $\nabla r_{min} =0$ and $
{\nabla}^{2} r_{min} \ge 0$. In view of (\ref{2.3}), we get
\begin{align*}
    H(r_{min})&=\frac{1}{r_{min}}\left[n-\frac{1}{r_{min}}{\nabla}^{2} r_{min}\right]
\end{align*}
Thus
\begin{align*}
    \frac{\partial }{\partial t}r_{min}&=-nf{r_{min}^{-1}} + f{r_{min}^{-2}}\nabla^{2}r_{min} - \frac{n}{n-1} \frac{\partial f}{\partial r}(r_{min})\\
    &\geq -nf{r_{min}^{-1}} - \frac{n}{n-1} \frac{\partial f}{\partial r}(r_{min})
\end{align*}
If $-nf{r_{min}^{-1}} - \frac{n}{n-1} \frac{\partial f}{\partial
r}(r_{min})\geq 0$, then
\begin{align*}
    \frac{\partial }{\partial t}r_{min}\geq 0
\end{align*}
and by the above differential inequality, we infer that
\begin{align*}
    r_{min}(t) \leq \left(\frac{f(r_{min}(0))}{f(r_{min}(t))}\right)^{\frac{1}{n-1}}r_{min}(0)
\end{align*}
and since $\frac{\partial f}{\partial r}(r_{min})\leq
-\frac{1}{n-1}fr_{min}^{-1}\leq 0$, we have
\begin{align*}
    \frac{f(r_{min}(0))}{f(r_{min}(t))} \geq 1
\end{align*}
Then
\begin{align*}
    r_{min}(t) \ge  r_{min}(0)
\end{align*}
If $-nf{r_{min}^{-1}} - \frac{n}{n-1} \frac{\partial f}{\partial
r}(r_{min})< 0$, then $\widehat{f}(r_{min})>0$. Since $\widehat{f}$
has a zero point, then there are constants $0<\varepsilon<\delta$
such that
\begin{align*}
    \widehat{f}(\varepsilon)<0;\qquad \widehat{f}(\delta)>0
\end{align*}
Thereby
\begin{align*}
    r_{min}(t) > \delta
\end{align*}
Combining the above two scenarios, we get
\begin{align*}
    r_{min}(t) \geq C_{1}:= max \{ \delta,r_{min}(0) \}
\end{align*}

In the same way, Set $r_{max}(t):=\underset{x\in
\mathbb{S}^{n}}{max}  r(\cdot,t)$, we have $\nabla r_{max} =0$ and $
{\nabla}^{2} r_{max} \leq 0$. Similarly, we obtain
\begin{align*}
    \frac{\partial r_{max}}{\partial t}&=-nfr_{max}^{-1} + fr_{max}^{-2} {\nabla}^{2}r_{max} - \frac{n}{n-1} \frac{\partial f}{\partial r}(r_{max})\\
    &\leq -nfr_{max}^{-1} - \frac{n}{n-1} \frac{\partial f}{\partial r}(r_{max})
\end{align*}
If $-nfr_{max}^{-1}-\frac{n}{n-1} \frac{\partial f}{\partial
r}(r_{max}) > 0$, we have $\widehat{f}(r_{max})<0$, thus
\begin{align*}
    r_{max}(t) < \varepsilon
\end{align*}
If $-nfr_{max}^{-1}-\frac{n}{n-1} \frac{\partial f}{\partial
r}(r_{max}) \leq 0$, then
\begin{align*}
    r_{max}(t) \leq min \left\{r_{max}(0),\left(\frac{f(r_{max}(0))}{f(r_{max}(t))}\right)^\frac{1}{n}r_{max}(0)\right\}
\end{align*}
When $\frac{\partial f}{\partial r}(r_{max})< 0$, we have
$\frac{f(r_{max}(0))}{f(r_{max}(t))} <  1$, then
\begin{align*}
    r_{max}(t) \leq \left(\frac{f(r_{max}(0))}{f(r_{max}(t))}\right)^\frac{1}{n}r_{max}(0):=\theta r_{max}(0),\quad  \theta< 1
\end{align*}
When $\frac{\partial f}{\partial r}(r_{max}) \geq 0$, we have
$\frac{f(r_{max}(0))}{f(r_{max}(t))} \geq  1$, then
\begin{align*}
    r_{max}(t) \leq  r_{max}(0)
\end{align*}
Hence
\begin{align*}
    r_{max}(t) \leq min \{r_{max}(0), \theta r_{max}(0)\} = \theta r_{max}(0)
\end{align*}
Combining the above two situations, we obtain
\begin{align*}
    r_{max}(t) \leq C_{2}:= min \{\varepsilon, \theta r_{max}(0)\}
\end{align*}
\hfill${\square}$

Next, we estimate the gradient of $r$ which is not very precise but
is sufficient.
\begin{lemma}\label{le4.2}
    Let $r \in {C^{\infty}(\mathbb{S}^{n} \times [0,T))}$ be a smooth solution to the initial value problem (\ref{4.1}). For any time $t \in [0,T)$, there is a positive constant $C$ depending on 1 and $M_{0}$, such that
    \begin{align}\label{4.3}
        \underset{x\in \mathbb{S}^{n}}{max} |\nabla r(\cdot,t)|  \leq C
    \end{align}
\end{lemma}
\noindent{\it \bf Proof.}~~Let $\omega=lnr$ and $\varphi =
\frac{1}{2}|\nabla \omega|^{2}$. We re-express geometric quantities
in terms of $\omega$.
\begin{align*}
    \begin{array}{l}
        g_{i j}=e^{2 \omega}\left(e_{i j}+\nabla_{i} \omega \nabla_{j} \omega\right) \\
        g^{i j}=e^{-2 \omega}\left(e_{i j}-\frac{\nabla_{i} \omega \nabla_{j} \omega}{1+|\nabla \omega|^{2}}\right), \\
        h_{i j}=e^{\omega}\left(1+|\nabla \omega|^{2}\right)^{-\frac{1}{2}}\left(e_{i j}+\nabla_{i} \omega \nabla_{i} \omega-\nabla_{i j} \omega\right) \\
        \end{array}
\end{align*}
and
\begin{align}\label{4.4}
    H=\sum g^{ij}h_{ij}= e^{- \omega} (1+|\nabla \omega|^{2})^{-\frac{1}{2}}\left[n-\left(e^{ij} - \frac{\nabla ^{i}\omega \nabla ^{j}\omega}{1+|\nabla \omega|^{2}}\right)\nabla_{ij}\omega\right]
\end{align}
From (\ref{4.1}), $\omega$ is the solution of the initial value
problem
\begin{align}\label{4.5}
    \begin{cases}
        \frac{\partial }{\partial t}\omega= {e^{-\omega}}\sqrt{1+|\nabla \omega|^{2}}\left(-fH-\frac{n}{n-1}e^{-\omega} \frac{\partial f}{\partial \omega}\sqrt{1+|\nabla \omega|^{2}}\right)\qquad \text{on} \quad \mathbb{S}^{n}\times[0,T)\\
        \omega(\cdot,0)=\omega_{0} = lnr_{0}
    \end{cases}
\end{align}
Combining $\frac{\partial}{\partial t}\varphi
=\nabla^{k}\omega\nabla_{k}\left(\frac{\partial }{\partial
t}\omega\right)$, we obtain the evolution equation of $\varphi$
\begin{align}
   \frac{\partial }{\partial t}\varphi=
   &-\nabla^{k}\omega \nabla_{k}(e^{-\omega})\sqrt{1+|\nabla \omega|^{2}}fH
   - e^{-\omega}{\nabla^{k}\omega} \nabla_{k}\left(\sqrt{1+|\nabla \omega|^{2}}\right)fH
   \notag \\
   & - e^{-\omega}\sqrt{1+|\nabla \omega|^{2}}{\nabla^{k}\omega}\nabla_{k}f H
   {-e^{- \omega}\sqrt{1+|\nabla \omega|^{2}} f{\nabla^{k}\omega} \nabla_{k}H}
   \notag\\
  &  -\frac{n}{n-1}{\nabla^{k}\omega} \nabla_{k}(e^{-2\omega})\frac{\partial f}{\partial \omega}(1+|\nabla \omega|^{2})
  -\frac{n}{n-1}e^{- 2\omega}{\nabla^{k}\omega}\nabla_{k}\left(\frac{\partial f}{\partial \omega}\right)(1+|\nabla \omega|^{2})
 \notag \\
  &  -\frac{n}{n-1}e^{-2 \omega}\frac{\partial f}{\partial \omega}{\nabla^{k}\omega} \nabla_{k}(1+|\nabla \omega|^{2})
   \label{4.6}
\end{align}
where
\begin{align*}
    \nabla_{k}f=\frac{\partial f}{\partial \omega}\nabla_{k}\omega; \qquad \nabla_{k}\left(\frac{\partial f}{\partial \omega}\right)=\frac{\partial^{2} f}{\partial \omega^{2}}\nabla_{k}\omega
\end{align*}
Suppose $\varphi$ attains the spatial maximum at point
$\left(x_{t},t\right) \in \left(\mathbb{S}^{n} \times [0,T)\right)$,
we have
\begin{align}
    \nabla \varphi&=\nabla^{m}\omega \nabla_{km}\omega =0  \label{4.7}\\
    \nabla^{2}\varphi&=\nabla^{km}\omega \nabla_{km}\omega + \nabla^{k}\omega \nabla_{klm}\omega  \leq 0\label{4.8} \qquad k=1,\dots ,n
\end{align}
and also
\begin{align*}
    \nabla_{k}\left(\sqrt{1+|\nabla \omega|^{2}}\right)=&\left(1+|\nabla \omega|^{2}\right)^{-\frac{1}{2}}\nabla^{m}\omega \nabla_{km}\omega =0\\
    H\nabla^{k}\omega  =& n e^{-\omega}\left(1+|\nabla \omega|^{2}\right)^{-\frac{1}{2}} \nabla^{k}\omega\\
        \nabla_{k}H=& -e^{\omega}\nabla_{k}\omega \left(1+|\nabla \omega|^{2}\right)^{-\frac{1}{2}} \left[n-\left(e^{ij} - \frac{\nabla ^{i}\omega \nabla ^{j}\omega}{1+|\nabla \omega|^{2}}\right)\nabla_{ij}\omega\right] \\
    &-e^{\omega}\left(e^{ij} - \frac{\nabla ^{i}\omega \nabla ^{j}\omega}{1+|\nabla \omega|^{2}}\right)\nabla_{kij}\omega\\
    \nabla^{k}\omega {\nabla_{k}H} \geq & -ne^{-\omega}|\nabla \omega|^{2} \left(1+|\nabla \omega|^{2}\right)^{-\frac{1}{2}}
\end{align*}
Bring (\ref{4.7}), (\ref{4.8}) and the above equations into
(\ref{4.6}), we obtain
\begin{align}
    \frac{\partial }{\partial t}\varphi &\leq e^{-2\omega}|\nabla \omega|^{2}\left[2nf-n\frac{\partial f}{\partial \omega}+\frac{2n}{n-1}\frac{\partial f}{\partial \omega}\left(1+|\nabla \omega|^{2}\right)-\frac{n}{n-1}\frac{\partial^{2} f}{\partial \omega^{2}}\left(1+|\nabla \omega|^{2}\right)\right] \label{4.9}
\end{align}
Since $\widehat{f}$ is monotonically increasing with respect to $r$
, we get
\begin{align}\label{4.10}
     2nf -n \frac{\partial f}{\partial r} r \leq \frac{n}{n-1}\frac{\partial^{2} f}{\partial r^{2}}r^{2} - \frac{n}{n-1} \frac{\partial f}{\partial r} r
\end{align}
Since
\begin{align*}
    \frac{\partial f}{\partial r}= r^{-1} \frac{\partial f}{\partial \omega};\qquad \frac{\partial^{2} f}{\partial r^{2}}=r^{-2}\left(\frac{\partial^{2} f}{\partial \omega^{2}}-\frac{\partial f}{\partial \omega}\right)
\end{align*}
then the inequality (\ref{4.10}) has the following form
\begin{align}\label{4.11}
     2nf -n \frac{\partial f}{\partial \omega}  \leq \frac{n}{n-1}\frac{\partial^{2} f}{\partial \omega^{2}}- \frac{2n}{n-1} \frac{\partial f}{\partial \omega}
\end{align}
If $\frac{n}{n-1}\frac{\partial^{2} f}{\partial \omega^{2}}-
\frac{2n}{n-1} \frac{\partial f}{\partial \omega} \geq 0$, we have
\begin{align*}
    2nf -n \frac{\partial f}{\partial \omega}  \leq \frac{n}{n-1}\frac{\partial^{2} f}{\partial \omega^{2}}-  \frac{2n}{n-1} \frac{\partial f}{\partial \omega} \leq  \left( \frac{n}{n-1}\frac{\partial^{2} f}{\partial \omega^{2}}- \frac{2n}{n-1} \frac{\partial f}{\partial \omega}\right) \left(1+|\nabla \omega|^{2}\right)
\end{align*}
Then
\begin{align*}
    \frac{\partial \varphi}{\partial t} \leq  0
\end{align*}
If $\frac{n}{n-1}\frac{\partial^{2} f}{\partial \omega^{2}}-
\frac{2n}{n-1} \frac{\partial f}{\partial \omega} < 0$, we have
\begin{align*}
    \frac{\partial \varphi}{\partial t} &\leq 2\varphi e^{-2 \omega}\left(2nf -n \frac{\partial f}{\partial \omega} \right)-2\varphi\left(1+2\varphi\right)e^{-2 \omega}\left(\frac{n}{n-1}\frac{\partial^{2} f}{\partial \omega^{2}} - \frac{2n}{n-1} \frac{\partial f}{\partial \omega} \right) \\
    &\leq2 \left(2nf + \frac{3n-n^{2}}{n-1} \frac{\partial f}{\partial \omega}  -\frac{n}{n-1}\frac{\partial^{2} f}{\partial \omega^{2}} \right)\varphi e^{-2 \omega}-4\left(\frac{n}{n-1}\frac{\partial^{2} f}{\partial \omega^{2} } - \frac{2n}{n-1} \frac{\partial f}{\partial \omega} \right)\varphi e^{-2 \omega*}
\end{align*}
Let
\begin{align*}
    C_{3}&:=2e^{-2 \omega}\left(2nf + \frac{3n-n^{2}}{n-1} \frac{\partial f}{\partial \omega}  -\frac{n}{n-1}\frac{\partial^{2} f}{\partial \omega^{2}} \right) \leq 0\\
    C_{4}&:=-4e^{-2 \omega}\left(\frac{n}{n-1}\frac{\partial^{2} f}{\partial \omega^{2}} - \frac{2n}{n-1} \frac{\partial f}{\partial \omega} \right)\geq 0
\end{align*}
then the inequality (\ref{4.9}) is expressed as
\begin{align}\label{4.12}
    \frac{\partial \varphi}{\partial t} \leq C_{3} \varphi + C_{4}\varphi^{2}
\end{align}
Assume that $ C_{3}+C_{4}\varphi (t) \geq 0 $, otherwise we have
$\varphi (t) < -\frac{C_{3}}{C_{4}}$ and $\frac{\partial }{\partial
t} \varphi <0$. Since
\begin{align*}
    -\frac{C_{3}}{C_{4}} = -\frac{2e^{-2 \omega}\left(2nf -n \frac{\partial f}{\partial \omega} +\frac{n}{n-1}\frac{\partial^{2} f}{\partial \omega^{2}} - \frac{2n}{n-1} \frac{\partial f}{\partial \omega} \right)}{-4e^{-2 \omega}\left(\frac{n}{n-1}\frac{\partial^{2} f}{\partial \omega^{2}} - \frac{n}{n-1} \frac{\partial f}{\partial \omega} \right)}
     = \frac{1}{2}\left(\frac{2nf-n\frac{\partial f}{\partial \omega}}{\frac{n}{n-1}\frac{\partial^{2} f}{\partial \omega^{2}}-\frac{n}{n-1}\frac{\partial f}{\partial \omega}}+1\right)\leq 1
\end{align*}
Hence
\begin{align*}
    \varphi < C:= min\{1,\varphi(0)\}
\end{align*}
We know that $\varphi (t) \geq -\frac{C_{3}}{C_{4}} $, This implies
there exists $t_{1} \in [0,T)$ such that
$\varphi_{1}:=\varphi(t_{1})=-\frac{C_{3}}{C_{4}}$, then
\begin{align*}
    \frac{\partial \varphi}{\partial t} <  C_{3}\varphi +\left(C_{4}-\frac{C_{4}}{C_{3}}\right)\varphi^{2}
\end{align*}
where $-\frac{C_{4}}{C_{3}} \geq 1$. Solving the above inequality
yields
\begin{align*}
    \frac{\varphi}{C_{3}\varphi +\left(C_{4}-\frac{C_{4}}{C_{3}}\right)} > -\frac{C_{3}}{C_{4}} e^{C_{3}t} \geq 0
\end{align*}
Thus
\begin{align*}
    \varphi < \frac{-C_{3}}{C_{4}-\frac{C_{4}}{C_{3}}-1} \leq -\frac{C_{3}}{C_{4}}
\end{align*}
This contradicts the hypothesis,  so $\varphi < C$ holds.
\hfill${\square}$

Finally, we establish the uniform positive upper bound of the
principal curvatures.
\begin{lemma}
    Along the flow (1.5), the principal curvatures of $M_{t}$ satisfy
    \begin{align}\label{4.14}
        \kappa_{i} \leq  \widetilde{C} ,\qquad i=1,\cdots,n
    \end{align}
    where $\widetilde{C}$ is a positive constant and depends only on $M_{0}$.
\end{lemma}
\noindent{\it \bf Proof.}~~Suppose that at time $t_{0} \in[0, T)$,
$\kappa$ attains its maximum at point $x_{0}$ in the direction
$\zeta_{0} \in T_{X(x_{0},t_{0})}M_{t_{0}}$. After rotating
coordinate, we may assume that $x_{0}\in \mathbb{S}^{n}$ in the
north pole and  $M_{t}$ can be represented as the graph of $u$ in
the neighborhood of $\left(X_{0},t_{0}\right) \subset \Sigma \times
[0,T)$, where $X_{0}=X(x_{0},t_{0})\subset M_{t_{0}}$, $\Sigma
\subset T_{X_{0}}M_{t_{0}}$ is the tangent hyperplane, and $f$ can
be expressed as $f(u)$. We know that
\begin{align}\label{4.15}
    \frac{\partial f}{\partial t} =\frac{\partial f}{\partial r} \frac{\partial r}{\partial t} =\frac{\partial f}{\partial u} D_{t}u
\end{align}
Bring (\ref{2.7}) and (\ref{2.9}) into (\ref{4.15}) respectively, we
can get
\begin{align*}
    \frac{\partial f}{\partial t}=\frac{\partial f}{\partial r}\left( \sqrt{ 1+ r^{-2}|\nabla r|^{2}}\right)\Phi(x,t)
\end{align*}
and
\begin{align*}
    \frac{\partial f}{\partial t}=\frac{\partial f}{\partial u} \left(-\sqrt{1+|Du|^{2}}\right)\Phi(x,t)
\end{align*}
Then, $\frac{\partial f}{\partial r}$ and $\frac{\partial
f}{\partial u}$ are related as follows
\begin{align}\label{4.16}
    \frac{\partial f}{\partial r}= - \frac{\partial f}{\partial u} \frac{\sqrt{1+|Du|^{2}}}{\sqrt{1+r^{-2}|\nabla r|^{2}}}
\end{align}
Therefore, $u$ is the solution of the following equation
\begin{align}\label{4.17}
    \begin{cases}
        D_{t}u=\sqrt{1+|Du|^{2}}\left( fH - \frac{n}{n-1}\frac{\partial f }{\partial u}\sqrt{1+|Du|^{2}}\right)\qquad \text{on} \quad \Sigma \times [0,T)\\
        u(\cdot,0)=u_{0}
    \end{cases}
\end{align}
For convenience, we choose an orthogonal coordinate system
$\{\eta_{1}, \ldots, \eta_{n}\}$ in $\Sigma$ with the origin at
$X_{0}$. Then in coordinates centered at the origin and parallel to
$\{\eta_{1}, \ldots, \eta_{n}\}$ we have
\begin{align*}
   X_{0}=\left(a_{1}, \ldots, a_{n}, a\right)
\end{align*}
for appropriate constants $a_{1}, \ldots, a_{n}, a$ and $a>0$. Thus
\begin{align*}
    &X_{t}=\left(a_{1}, \ldots, a_{n}, a\right)+(\eta, u(\eta)),
    \\
    &\nu =\frac{(D u,-1)}{v}
\end{align*}
where $v=\left(1+|D u|^{2}\right)^{\frac{1}{2}}$. Without loss of
generality, we can assume that the maximum principal curvature of
the graph $u(\cdot,t)$ is obtained at $(0,t_{0})$ in the direction
$\eta_{1}$. In the vicinity of $(0,t_{0})$, we have
\begin{align}\label{4.18}
    \kappa_{1}=\frac{D_{11} u}{v\left(1+\left|D_{1} u\right|^{2}\right)} .
\end{align}
and
\begin{align}
    u(0, t_{0})=|D u(0, t_{0})|&=0 \label{4.19}\\
    D_{1 \alpha} u(0, t_{0})&=0 \qquad \text { for } \alpha>1 \label{4.20}
\end{align}
We can rotate $\{\eta_{2}, \ldots, \eta_{n}\}$ such that  $D^{2}
u(0, t_{0})$ is diagonal and $D_{11} u(0, t_{0})>0$. In the
equations above and below, the convention of summing over Latin
indices will be used instead of Greek indices. By calculating
\begin{align}\label{4.21}
    D_{t}\kappa_{1}=
    \frac{D_{11t}u}{v\left(1+\left|D_{1} u\right|^{2}\right)}
    -\frac{2D_{1}uD_{1t}uD_{11}u}{v \left(1+\left|D_{1} u\right|^{2}\right)^{2}}
    -\frac{v^{-1}D^{k}uD_{kt}u}{v^{2}\left(1+\left|D_{1} u\right|^{2}\right)}
\end{align}
Then, at $(0,t_{0})$
\begin{align*}
    D_{t}\kappa_{1}=D_{11t}u
\end{align*}
We also have
\begin{align*}
    D_{\alpha}\kappa_{1}=
    \frac{D_{11\alpha}u}{v\left(1+\left|D_{1} u\right|^{2}\right)}
    -\frac{2D_{1}uD_{1\alpha}uD_{11}u}{v \left(1+\left|D_{1} u\right|^{2}\right)^{2}}
    -\frac{v^{-1}D^{k}u D_{k\alpha}u}{v^{2}\left(1+\left|D_{1} u\right|^{2}\right)}
\end{align*}
and
\begin{align*}
    D_{\alpha \alpha}\kappa_{1}= &\frac{ D_{11\alpha \alpha}u}{v\left(1+\left|D_{1} u\right|^{2}\right)}
   -\frac{D_{11\alpha}u D_{\alpha}v)}{v^{2}\left(1+\left|D_{1} u\right|^{2}\right)}
   -\frac{2D_{11\alpha}u D^{1}uD_{1\alpha}u}{v\left(1+\left|D_{1} u\right|^{2}\right)^{2}} \notag\\
   &-\frac{D_{11\alpha}u D_{\alpha}v+D_{11}u D_{\alpha \alpha}v}{v^{2}\left(1+\left|D_{1} u\right|^{2}\right)} \notag \\
   &+\frac{D_{11}u D_{\alpha}v}{v^{4}\left(1+\left|D_{1} u\right|^{2}\right)^{2}}\left[2vD_{\alpha}v \left(1+\left|D_{1} u\right|^{2}\right) + 2v^{2}D^{1}uD_{1\alpha}u\right] \notag \\
   &+\frac{2D_{11}u D^{1}uD_{1\alpha}u}{v^{2}\left(1+\left|D_{1} u\right|^{4}\right)}
   \left[\left(1+\left|D_{1} u \right|^{2}\right) D_{\alpha}v + 4v\left(1+\left|D_{1} u\right|^{2}\right)D^{1}uD_{1\alpha}u \right]  \notag\\
   &-\frac{2}{v \left(1+\left|D_{1} u\right|^{2}\right)^{2}}\left[ D_{11\alpha}u D^{1}u D_{1\alpha}u
   + D_{11}u D^{1\alpha}u D_{1\alpha}u  + D_{11}u D^{1}uD_{1\alpha \alpha}u\right]
\end{align*}
where
\begin{align*}
    D_{\alpha}v&=v^{-1}D^{k}uD_{k\alpha}u \\
    D_{\alpha \alpha}v&=-v^{-3}\left(D^{k}uD_{k\alpha}u\right)^{2}+v^{-1}\left(D_{\alpha \alpha}u\right)^{2}+v^{-1}D^{k}uD_{k\alpha\alpha}u
\end{align*}
Thus, at $(0,t_{0})$
\begin{align}
    D_{\alpha} v&=0 \label{4.22}\\
    D_{\alpha \alpha} v&=\left(D_{\alpha \alpha} u\right)^{2} \label{4.23}
\end{align}
and
\begin{align}
    D_{\alpha}\kappa_{1} &= D_{11\alpha}u =0 \label{4.24}\\
    D_{\alpha \alpha}\kappa_{1}&=D_{\alpha \alpha}(D_{11}u)- D_{11}u( D_{\alpha \alpha}v)-2D_{11}u(D|D_{1}u|)^{2} \notag \\
    &=D_{11\alpha \alpha}u - D_{11}u( D_{\alpha \alpha}u)^{2}- 2(D_{11}u)^{3}\leq 0 \label{4.25}
\end{align}
Next, we compute $ D_{11t}u$. Differentiating (\ref{4.18}) twice in
$\eta_{1}$ direction, we get
\begin{align*}
    D_{11t}u = & D_{11}v \left( fH - \frac{n}{n-1}\frac{\partial f }{\partial u}v\right) +2D_{1}v\left(HD_{1}f+fD_{1}H  - \frac{n}{n-1}\frac{\partial f}{\partial u}D_{1}v\right)\\
    &+v\left[HD_{11}f+fD_{11}H + 2D_{1}f D_{1}H - \frac{n}{n-1}vD_{11}\left(\frac{\partial f}{\partial u}\right)  - \frac{n}{n-1}\frac{\partial f}{\partial u}D_{11}v \right]\\
    &-4 \frac{n}{n-1} v D_{1}\left(\frac{\partial f}{\partial u}\right)
\end{align*}
and at $(0,t_{0})$
\begin{align}
    D_{11t}u= &(D_{11}u)^{2}\left( fH - \frac{n}{n-1}\frac{\partial f }{\partial u}\right) -\frac{n}{n-1}D_{11}\left(\frac{\partial f}{\partial u}\right)- \frac{n}{n-1}\frac{\partial f}{\partial u}(D_{11}u)^{2} \notag \\
    &+HD_{11}f+ fD_{11}H +2D_{1}fD_{1}H \label{4.26}
\end{align}
From $H=S_{1}(\kappa )=\kappa_{1}+ \cdots +\kappa_{n}$, we compute
that
\begin{align*}
  D_{1}H&=S_{1}^{ij}D_{1}\lambda_{ij}=D_{1}\lambda_{ij}\\
  D_{11}H&= S_{1}^{ij,pq}D_{1}\lambda_{ij}D_{1}\lambda_{pq} + S_{1}^{ij}D_{11}\lambda_{ij}=S_{1}^{ij}D_{11}\lambda_{ij}
\end{align*}
at $(0,t_{0})$. By direct calculation,
\begin{align*}
    D_{1}f&=\frac{\partial f}{\partial u}D_{1}u=0 \qquad D_{11}f=\frac{\partial^{2}f}{\partial u^{2}}(D_{1}u)^{2} + \frac{\partial f}{\partial u}D_{11}u =\frac{\partial f}{\partial u}D_{11}u
\end{align*}
and
\begin{align*}
    D_{1}(\frac{\partial f}{\partial u}) = \frac{\partial^{2}f}{\partial u^{2}}D_{1}u =0
    \qquad D_{11}(\frac{\partial f}{\partial u}) = \frac{\partial^{3}f}{\partial u^{3}}(D_{1}u)^{2} +\frac{\partial^{2}f}{\partial u^{2}}D_{11}u  = \frac{\partial^{2}f}{\partial u^{2}}D_{11}u
\end{align*}
at $(0,t_{0})$. Bring the above formulas into (\ref{4.26}), we get
\begin{align}
    D_{11t}u=&(D_{11}u)^{2}\left( fH - \frac{n}{n-1}\frac{\partial f }{\partial u}\right) + \frac{\partial f}{\partial u}HD_{11}u + fS_{1}^{ij}D_{11}\lambda_{ij} \notag \\
    &-\frac{n}{n-1}\frac{\partial^{2} f}{\partial u^{2}}(D_{11}u)- \frac{n}{n-1}\frac{\partial f}{\partial u}(D_{11}u)^{2} \label{4.27}
\end{align}
and the equation (\ref{2.8}) implies that
\begin{align}\label{4.28}
    D_{11}\lambda_{ij}=D_{11 i j} u-D_{1 k} u D_{1 k} u D_{i j} u-2 D_{1 i} u D_{1 k} u D_{j k} u .
\end{align}
at $(0,t_{0})$. From $ D^{2} u$ and $\left[S_{1}^{ij}\right]$ are
diagonal at $(0,t_{0})$, we have
\begin{align*}
    &D_{11} \lambda_{11}=D_{1111} u-3\left(D_{11} u\right)^{3}\\
    &D_{11} \lambda_{ii}=D_{11 ii} u-D_{ii} u\left(D_{11} u\right)^{2}\\
    &S_{1}^{ii}D_{11}\lambda_{ii}=D_{11}\lambda_{11}+\cdots+D_{nn}\lambda_{nn}
\end{align*}
So by (\ref{4.25}),
\begin{align}\label{4.29}
    D_{11ii}u \leq D_{11}u( D_{ii}u)^{2}+ 2(D_{11}u)^{3}
\end{align}
Taking (\ref{4.29}) into (\ref{4.28}), we get
\begin{align*}
    D_{11}\lambda_{ii} \leq D_{11}u( D_{ii}u)^{2}+ 2(D_{11}u)^{3} -D_{ii}u(D_{11}u)^{2}
\end{align*}
Then
\begin{align}
    S_{1}^{ii}D_{11}\lambda_{ii} &\leq D_{11}u S_{1}^{ii}\left((D_{ii}u)^{2}
    + 2(D_{11}u)^{2}-D_{ii}u D_{11}u\right) \notag \\
    &\leq (n+2)(D_{11}u)^{3} -H(D_{11}u)^{2} \label{4.30}
\end{align}
Using (\ref{4.27}) and (\ref{4.30}) in (\ref{4.21}), we obtain, at
$(0,t_{0})$
\begin{align}\label{4.31}
    D_{t}\kappa_{1} \leq  &(n+2)f(D_{11}u)^{3}+\frac{\partial f}{\partial u}HD_{11}u-\frac{n}{n-1}\frac{\partial^{2} f}{\partial u^{2}}D_{11}u- \frac{2n}{n-1}\frac{\partial f}{\partial u}(D_{11}u)^{2}
\end{align}
Now extended the condition (\ref{4.10}) to $f[u(\eta,t)]$. According
to the different representations of $f$
\begin{align*}
    &f:=f(r)=f\circ r \circ \xi :M_{t} \to \mathbb{S}^{n} \to \mathbb{R}; \quad
    \xi:M_{t} \to \mathbb{S}^{n};\quad
    X_{t}(x)=r(\xi,t)\xi\\
    &f:=f(u)=f\circ u \circ \eta :M_{t} \to \Sigma \to \mathbb{R};\quad
    \eta :M_{t} \to \Sigma;\quad X_{t}(x) =(\eta,u(\eta,t))
\end{align*}
we have
\begin{align*}
    \frac{\partial f}{\partial x} =\frac{\partial f}{\partial r}\frac{\partial r}{\partial \xi}\frac{\partial \xi}{\partial x} =\frac{\partial f}{\partial r}\nabla r\frac{\partial \xi}{\partial x}
\end{align*}
and
\begin{align*}
    \frac{\partial f}{\partial x}=\frac{\partial f}{\partial u}\frac{\partial u}{\partial \eta}\frac{\partial \eta}{\partial x} =\frac{\partial f}{\partial u}Du\frac{\partial \eta}{\partial x}
\end{align*}
By (\ref{4.19}) we have $|\frac{\partial f}{\partial x}| =0$, i.e.
$\left|\frac{\partial f}{\partial r}\nabla r\frac{\partial
\xi}{\partial x}\right|=0$ at  $(0, t_{0})$. There is either
$\frac{\partial f}{\partial r}= 0$ or $\nabla r\frac{\partial
\xi}{\partial x}=0$.

Case 1: If $\frac{\partial f}{\partial r} = 0$, from (\ref{4.18}),
we obtain $\frac{\partial f}{\partial u} = 0$. Then, at $(0,t_{0})$,
\begin{align}\label{4.32}
    D_{t}\kappa_{1}\leq (n+2)f(D_{11}u)^{3} -\frac{n}{n-1}\frac{\partial^{2}f}{\partial u^{2}}D_{11}u
\end{align}
and the inequality (\ref{4.10}) becomes
\begin{align}\label{4.33}
    0 &\leq 2nf \leq \frac{n}{n-1}r_{0}^{2}\frac{\partial^{2}f}{\partial r^{2}}
\end{align}
where $r_{0}:=r^{2}(0,t_{0}) = a_{1}^{2}+\cdots +a_{n}^{2} +a^{2}$.
Similarly, we have
\begin{align*}
    \frac{\partial^{2}f}{\partial t^{2}} =\frac{\partial^{2}f}{\partial r^{2}}\left(\frac{\partial r}{\partial t}\right)^{2}+ \frac{\partial f}{\partial r} \frac{\partial^{2}r}{\partial t^{2}}= \frac{\partial^{2}f}{\partial u^{2}}\left(\frac{\partial u}{\partial t}\right)^{2}+ \frac{\partial f}{\partial u} \frac{\partial^{2}u}{\partial t^{2}}
\end{align*}
and
\begin{align}\label{4.34}
    \frac{\partial^{2}f}{\partial t^{2}} &= \frac{\partial^{2}f}{\partial r^{2}}\left(\frac{\partial r}{\partial t}\right)^{2}=\frac{\partial^{2}f}{\partial u^{2}}\left(\frac{\partial u}{\partial t}\right)^{2}
\end{align}
at $(0,t_{0})$. Using (\ref{4.15}) and (\ref{4.16}), we get
\begin{align}\label{4.35}
    \frac{\partial u}{\partial t}= -\left(\frac{\sqrt{1+|Du|^{2}}}{\sqrt{1+r^{-2}|\nabla r|^{2}}}\right)\frac{\partial r}{\partial t}
\end{align}
Thus, at $(0,t_{0})$,
\begin{align*}
    \frac{\partial^{2}f}{\partial r^{2}}=\frac{\partial^{2}f}{\partial u^{2}}\frac{1}{1+r^{-2}|\nabla r|^{2}} \leq \frac{\partial^{2}f}{\partial u^{2}}
\end{align*}
and $\frac{\partial^{2}f}{\partial u^{2}}\geq 0$ can be obtained
from (\ref{4.33}). Hence, $f(u)$ satisfies the following condition
\begin{align}\label{4.36}
    0\leq 2nf \leq \frac{n}{n-1}r_{0}^{2} \frac{\partial^{2}f}{\partial u^{2}}
\end{align}
at $(0,t_{0})$. Bring (\ref{4.36}) into (\ref{4.32}), we have, at
$(0,t_{0})$,
\begin{align*}
    \frac{\partial \kappa_{1}}{\partial t}
    &\leq \frac{1}{2} \frac{n+2}{n-1}  \frac{\partial^{2}f}{\partial u^{2}} r_{0}^{2}(D_{11}u)^{3} -\frac{n}{n-1}\frac{\partial^{2}f}{\partial u^{2}}D_{11}u \\
    &=D_{11}u\frac{\partial^{2}f}{\partial u^{2}}\left( \frac{1}{2} \frac{n+2}{n-1} r_{0}^{2}(D_{11}u)^{2} - \frac{n}{n-1}\right)
\end{align*}
We can choose $X_{0}$ to be close enough to the original origin such
that $r_{0}$ is small enough and $r_{0}^{2}(D_{11}u)^{2} \leq
\epsilon$, where $\frac{1}{4} < \epsilon < \frac{2n}{n+2}$. Thus, we
get
\begin{align*}
    D_{t}\kappa_{1} \leq 0
\end{align*}

Case 2: If $\nabla r\frac{\partial \xi}{\partial x}=0$, we obtain,
at $(0,t_{0})$
\begin{align*}
    \nabla r\frac{\partial \xi}{\partial x}=\frac{\partial r}{\partial x} =\frac{\partial |X(x,t)|}{\partial x}  =|X|^{-1}\left<\frac{\partial X}{\partial x} ,X\right> =\left<\frac{\partial X}{\partial x} ,\xi \right> =0
\end{align*}
From the above equation, we deduce that the unit normal vector is
proportional to $\xi$, so we have
\begin{align}\label{4.37}
    \nabla r= 0
\end{align}
at $(0,t_{0})$. Combining with (\ref{4.16}) and (\ref{4.35}), we
also have
\begin{align}\label{4.38}
    \frac{\partial f}{\partial u}= - \frac{\partial f}{\partial r};\qquad
    \frac{\partial u}{\partial t}= - \frac{\partial r}{\partial t}
\end{align}
Through direct calculation, we can get
\begin{align*}
    \frac{\partial^{2}r}{\partial t^{2}} =& \frac{\partial }{\partial t}\left(- \frac{\sqrt{1+r^{-2}|\nabla r|}^{2}}{\sqrt{1+|Du|^{2}}} \frac{\partial u}{\partial t}\right) \\
    =&-\left(1+|Du|^{2}\right)^{-1}\sqrt{1+|Du|^{2}}\frac{\partial}{\partial t}\left(\sqrt{1+r^{-2}|\nabla r|^{2}}\right)\frac{\partial u}{\partial t}\\
    &-\sqrt{1+r^{-2}|\nabla r|^{2}}\frac{\partial}{\partial t}\left(\sqrt{1+|Du|^{2}}\right)\frac{\partial u}{\partial t} \\
    &-\left(\frac{\sqrt{1+r^{-2}|\nabla r|^{2}}}{\sqrt{1+|Du|^{2}}}\right)\frac{\partial^{2}u}{\partial t^{2}}
\end{align*}
where
\begin{align*}
    \frac{\partial}{\partial t}\left(\sqrt{1+|Du|^{2}}\right)&=\left(1+|Du|^{2}\right)^{-\frac{1}{2}}DuD_{t}(Du) \\
    \frac{\partial}{\partial t}\left(\sqrt{1+r^{-2}|\nabla r|^{2}}\right)&=\left(1+r^{-2}|\nabla r|^{2}\right)^{-\frac{1}{2}} r^{-4}\left(r^{2}\nabla r\frac{\partial}{\partial t}|\nabla r-r|\nabla r|^{2}\frac{\partial }{\partial t}r\right)
\end{align*}
and
\begin{align*}
    \frac{\partial}{\partial t}\left(\sqrt{1+|Du|^{2}}\right)=0\\
    \frac{\partial}{\partial t}\left(\sqrt{1+r^{-2}|\nabla r|^{2}}\right)=0
\end{align*}
at $(0,t_{0})$. Then
\begin{align}\label{4.39}
    \frac{\partial^{2}r}{\partial t^{2}} =-\frac{\sqrt{1+r^{-2}|\nabla r|^{2}}}{\sqrt{1+|Du|^{2}}}\frac{\partial^{2}u}{\partial t^{2}}= -\frac{\partial^{2}u}{\partial t^{2}}
\end{align}
Taking (\ref{4.38}) and (\ref{4.39}) into (\ref{4.34}), we get, at
$(0,t_{0})$,
\begin{align}\label{4.40}
    \frac{\partial^{2}f}{\partial r^{2}} =\frac{\partial^{2}f}{\partial u^{2}}
\end{align}
Therefore, combining (\ref{4.10}) and the above relation equations,
we have
\begin{align}\label{4.41}
    2nf + nr_{0}\frac{\partial f}{\partial u} \leq \frac{n}{n-1}r_{0}^{2}\frac{\partial^{2}f}{\partial u^{2}} +
    \frac{n}{n-1}r_{0}\frac{\partial f}{\partial u}
\end{align}
at $(0,t_{0})$. Bring (\ref{4.41}) into (\ref{4.31}) we get
\begin{align*}
    D_{t} \kappa_{1} &\leq \left(\frac{1}{2} \frac{n+2}{n-1}r_{0}^{2}\frac{\partial^{2}f}{\partial u^{2}}+\frac{1}{2} \frac{n+2}{n-1}(2-n)r_{0}\frac{\partial f}{\partial u}\right)(D_{11}u)^{3}\\
    &+HD_{11}u\frac{\partial f}{\partial u} -\frac{n}{n-1}\frac{\partial^{2}f}{\partial u^{2}}D_{11}u-\frac{2n}{n-1}\frac{\partial f}{\partial u}(D_{11}u)^{2}
\end{align*}
It is know that the following formula is valid
\begin{align*}
    \frac{\partial^{2}f}{\partial x^{2}}&=\frac{\partial^{2}f}{\partial r^{2}}\left(\frac{\partial r}{\partial \xi}\frac{\partial \xi}{\partial x}\right)^{2}+\frac{\partial f}{\partial r} \frac{\partial^{2}r}{\partial \xi^{2}}\left(\frac{\partial \xi}{\partial x}\right)^{2}+\frac{\partial f}{\partial r}\frac{\partial r}{\partial \xi}\left(\frac{\partial^{2}\xi}{\partial x^{2}}\right)\\
    &=\frac{\partial^{2}f}{\partial u^{2}}\left(\frac{\partial u}{\partial \eta}\frac{\partial \eta}{\partial x}\right)^{2}+\frac{\partial f}{\partial u} \frac{\partial^{2}u}{\partial \eta^{2}}\left(\frac{\partial \eta}{\partial x}\right)^{2}+\frac{\partial f}{\partial u}\frac{\partial u}{\partial \eta}\left(\frac{\partial^{2}\eta}{\partial x^{2}}\right)
\end{align*}
and the equation (\ref{4.19}) and (\ref{4.37}) implies
\begin{align*}
    \frac{\partial^{2}f}{\partial x^{2}}=\frac{\partial f}{\partial r} \frac{\partial^{2}r}{\partial \xi^{2}}\left(\frac{\partial \xi}{\partial x}\right)^{2}=\frac{\partial f}{\partial u} \frac{\partial^{2}u}{\partial \eta^{2}}\left(\frac{\partial \eta}{\partial x}\right)^{2}
\end{align*}
at $(0,t_{0})$. Thus
\begin{align*}
    \frac{\Delta r}{D^{2}u}=\frac{\frac{\partial^{2}r}{\partial \xi^{2}}}{\frac{\partial^{2}u}{\partial \eta^{2}}} =\frac{\frac{\partial f}{\partial u} \left(\frac{\partial \eta}{\partial x}\right)^{2}}{\frac{\partial f}{\partial r}\left(\frac{\partial \xi}{\partial x}\right)^{2}}
    =-\frac{\left(\frac{\partial \eta}{\partial x}\right)^{2}}{\left(\frac{\partial \xi}{\partial x}\right)^{2}}
    =-\left(\frac{\partial \eta}{\partial \xi}\right)^{2}
\end{align*}
From the definition of $\xi$, it follows that $\xi = \frac{X}{|X|} =
\frac{(\eta ,u)}{\sqrt{\eta ^{2}+u^{2}}}$. Differentiating $\xi$
with respect to $\eta$
\begin{align*}
    \frac{\partial \xi}{\partial \eta} =\frac{(1 ,Du)}{\sqrt{\eta ^{2}+u^{2}}} - \frac{(\eta,u)}{\left(\eta ^{2}+u^{2} \right)^{\frac{3}{2}}}(\eta+uDu)
\end{align*}
and
\begin{align*}
    \frac{\partial \xi}{\partial \eta} = \frac{a}{r_{0}^{3}}(a,-a_{1}, \cdots, -a_{n})
\end{align*}
at $(0,t_{0})$. Therefore
\begin{align}\label{4.42}
    \frac{\Delta r}{D^{2}u}=\frac{\frac{\partial^{2}r}{\partial \xi^{2}}}{\frac{\partial^{2}u}{\partial \eta^{2}}}=-\left(\frac{\partial \eta}{\partial \xi}\right)^{2} =-\left(\frac{\partial \xi}{\partial \eta}\right)^{-2} =-\frac{r_{0}^{4}}{a^{2}}
\end{align}
Using  $r=|X|$, we calculate that
\begin{align*}
    \Delta r&=-|X|^{2}\left<\nabla X,X\right>\left<\nabla X,X\right> + |X|^{-1}\left<\nabla^{2} X,X\right> + |X|^{-1}\left<\nabla X,\nabla X\right>\\
    &=-(\nabla r)^{2} + |X|^{-1}\left<\nabla^{2} X,X\right> + |X|^{-1}\left<\nabla X,\nabla X\right>
\end{align*}
Thus
\begin{align*}
    \Delta r  =|X|^{-1}\left<\nabla^{2} X,X\right> + |X|^{-1}\left<\nabla X,\nabla X\right>
\end{align*}
at $(0,t_{0})$. From $X=(\eta ,u(\eta))$, we have
\begin{align*}
    DX&=\frac{\partial X}{\partial \eta} =\nabla X \frac{\partial \xi}{\partial \eta}\\
    D^{2}X&=\nabla^{2}X \left(\frac{\partial \xi}{\partial \eta}\right)^{2}+\nabla X \frac{\partial^{2} \xi}{\partial \eta^{2}}\\
    \left< DX, DX\right>&=\left< \nabla X, \nabla X\right>\left(\frac{\partial \xi}{\partial \eta}\right)^{2}\\
    \left< D^{2}X, X\right>&=\left< \nabla^{2}X, X\right>\left(\frac{\partial \xi}{\partial \eta}\right)^{2} + \left< \nabla X, X\right>\frac{\partial^{2} \xi}{\partial \eta^{2}}
\end{align*}
and
\begin{align*}
    \left< D^{2}X, X\right>=\left< \nabla^{2}X, X\right>\left(\frac{\partial \xi}{\partial \eta}\right)^{2}
\end{align*}
at $(0,t_{0})$. Hence
\begin{align}\label{4.43}
  \Delta r= r_{0}^{-1}\left( \frac{\left< D^{2}X, X\right>}{\left(\frac{\partial \xi}{\partial \eta}\right)^{2}}+\frac{\langle DX, DX\rangle}{\left(\frac{\partial \xi}{\partial \eta}\right)^{2}}\right) =r_{0}^{-1}\left(\frac{\partial \eta}{\partial \xi}\right)^{2}(-aD^{2}u+1)
\end{align}
By (\ref{4.42}) and (\ref{4.43}), we get
\begin{align*}
    \Delta r = \frac{r_{0}^{4}}{a^{2}(r_{0}+a)}
\end{align*}
at $(0,t_{0})$. Thus
\begin{align*}
    H(0,t_{0}) = \frac{1}{r_{0}} \left(n-\frac{r_{0}^{3}}{a^{2}(r_{0}+a)}\right)
\end{align*}
We may assume that $a^{2} = \frac{1}{c^{2}}r_{0}^{2}$ where $c \in
\left(1,\frac{4}{3}\right)$, then
\begin{align*}
  0<(n-1)r_{0}^{-1} < H(0,t_{0}) < nr_{0}^{-1}
\end{align*}
By selecting a suitable $X_{0}$ can make $ \frac{1}{4} <
\left(r_{0}D_{11}u\right)^{2}  < \epsilon $, then
\begin{align*}
    D_{t} \kappa_{1}
    \leq &\frac{1}{2}\frac{n+2}{n-1}\left(\frac{\partial^{2}f}{\partial u^{2}}D_{11}u - (n-2)\frac{\partial f}{\partial u}r_{0}^{-1}D_{11}u\right)r_{0}^{2}(D_{11}u)^{2}
    +\frac{\partial f}{\partial u}D_{11}uH \\
    &- \frac{n}{n-1}\frac{\partial^{2}f}{\partial u^{2}}D_{11}u
    - 2\frac{n}{n-1}\frac{\partial f}{\partial u}(D_{11}u)^{2}\\
    \leq& - \frac{n}{n-1}(n-2)\frac{\partial f}{\partial u}r_{0}^{-1}D_{11}u
    +\frac{\partial f}{\partial u}D_{11}uH
    -2\frac{n}{n-1}\frac{\partial f}{\partial u}(D_{11}u)^{2}\\
    \leq& \frac{\partial f}{\partial u} \frac{D_{11}u}{r_{0}} \left( \frac{n}{n-1}
    - \frac{2n}{n-1}r_{0}D_{11}u\right)\leq 0
\end{align*}
Furthermore, we can apply the maximum principle to conclude that
$\kappa$ is bounded from above. \hfill${\square}$

Combining the above Lemmas we can get the longtime existence of flow
(\ref{1.5}).

\begin{proposition}
    Let $T^{*}$ be the maximal existence time of the flow (\ref{1.5}). Then $T^{*}=\infty$.
\end{proposition}
\noindent{\it \bf Proof.}~~From (\ref{4.1}), we know that
\begin{align*}
    \frac{\partial}{\partial t}r=-\sqrt{1+r^{-2}|\nabla r|^{2}}fH-\frac{n}{n-1} \frac{\partial f}{\partial r}(1+r^{-2}|\nabla r|^{2}):=\widetilde{\Phi}
\end{align*}
where $\widetilde{\Phi}:=\widetilde{\Phi}(\nabla r,\nabla^{2}r)$,
and $Q^{ij}=\frac{\partial \widetilde{\Phi}}{\partial
\left(\nabla_{ij}r\right)}$ is positive. Hence, (\ref{1.6}) is a
parabolic equation on $\mathbb{S}^{n} \times \mathbb{R}_{+}$, the
short time existence of the flow (\ref{1.5}) can be derived from the
theory of parabolic equation. By Proposition 4.1 and Proposition
4.2, we have the $C^{0}$ and $C^{1}$ estimate for the flow
(\ref{1.5}). In addition, from Proposition 4.3 we get the upper
bounded of the principal curvaturse, thus we obtain the uniform
$C^{2}$ estimate of the flow (\ref{1.5}). After that, applying
Krylov's \cite{K} theory and standard parabolic Schauder estimate to
derive $C^{2,\alpha}$ estimate, and higher order regular estimates
respectively. Therefore, we get the longtime existence of the flow
(\ref{1.5}). \hfill${\square}$

Finally, we prove the flow (\ref{1.5}) converges to the unique
sphere by making an exact estimate of the gradient of $r$.
\begin{proposition}\label{prop4.5}
    Let $r \in {C^{\infty}(\mathbb{S}^{n} \times [0,\infty))}$ be a smooth solution to the initial value problem (\ref{4.1}). For any time $t \in [0,\infty)$, there is a positive constant $\overline{C} $ depending only on $M_{0}$, such that
    \begin{align}\label{4.44}
        \underset{\xi\in \mathbb{S}_{n}}{max} |\nabla r(\cdot, t)|\leq \overline{C}e^{-\gamma t}
    \end{align}
\end{proposition}
\noindent{\it \bf Proof.}~~From Lemma \ref{le4.2}
\begin{align*}
    \frac{\partial }{\partial t}\varphi &\leq e^{-2\omega}|\nabla \omega|^{2}\left[2nf-n\frac{\partial f}{\partial \omega}+\frac{2n}{n-1}\frac{\partial f}{\partial \omega}\left(1+|\nabla \omega|^{2}\right)-\frac{n}{n-1}\frac{\partial^{2} f}{\partial \omega^{2}}\left(1+|\nabla \omega|^{2}\right)\right]
\end{align*}

If $\frac{n}{n-1}\frac{\partial^{2} f}{\partial \omega^{2}}-
\frac{2n}{n-1} \frac{\partial f}{\partial \omega} \geq 0$, we have
\begin{align*}
    \frac{\partial \varphi}{\partial t}
    \leq e^{-2\omega}|\nabla \omega|^{2}\left(2nf + \frac{3n-n^{2}}{n-1} \frac{\partial f}{\partial \omega}  -\frac{n}{n-1}\frac{\partial^{2} f}{\partial \omega^{2}} \right)
\end{align*}
Let
\begin{align*}
    C_{6}:&=e^{-2\omega}\left(2nf + \frac{3n-n^{2}}{n-1} \frac{\partial f}{\partial \omega}  -\frac{n}{n-1}\frac{\partial^{2} f}{\partial \omega^{2}} \right)\\
    &=r^{-2}\left(2nf  +\frac{2n-n^{2}}{n-1}\frac{\partial f}{\partial r} r - \frac{n}{n-1} \frac{\partial ^{2}f}{\partial r^{2}} r^{2}\right)
    \leq 0
\end{align*}
and since $f \in C^{\infty}(M)$, so $\frac{\partial f}{\partial r}$
and $\frac{\partial^{2} f}{\partial r^{2}}$ are continuous
functions. From Proposition 4.1, we know that $r$ is bounded, then
$\frac{\partial f}{\partial r}$ and $\frac{\partial^{2} f}{\partial
r^{2}}$ are bounded. Therefore, $C_{6}$ is bounded.

For some positive constant $\gamma \geq -\frac{1}{2} maxC_{6}$, we
obtain
\begin{align*}
    \frac{\partial \varphi}{\partial t} \leq -\gamma \varphi
\end{align*}
and this implies that the gradient of the radial function decreases
exponentially.

If $\frac{n}{n-1}\frac{\partial^{2} f}{\partial \omega^{2}}-
\frac{2n}{n-1} \frac{\partial f}{\partial \omega} < 0$, we have
\begin{align*}
    \frac{\partial \varphi }{\partial t} \leq C_{3}\varphi + C_{4} \varphi^{2}
\end{align*}
where $C_{3}+C_{4}\varphi<0$. There exists a constant $\gamma_{1}>0$
such that $C_{3}+C_{4}\varphi \leq -\gamma_{1}$, then
\begin{align*}
    \frac{\partial \varphi }{\partial t} \leq -\gamma_{1}\varphi
\end{align*}
So far the proof of Proposition \ref{prop4.5} is complete.
\hfill${\square}$

It can be inferred from Proposition \ref{prop4.5} that
$r_{\infty}:=\lim\limits_{t \to \infty} r(\cdot,t) $ exists and
\begin{align*}
    \underset{\xi\in \mathbb{S}_{n}}{max} |r(\cdot, t)-r_{\infty}|\leq Ce^{-\gamma t}
\end{align*}
Now we can apply the interpolation inequality and Sobolev embedding
theorem on $\mathbb{S}^{n}$ to derive the convergence of $r(\cdot,
t)$ to $r_{\infty}$ is smooth. From the comparison principle, the
uniquess of $r_{\infty}$ follows in a standard way. This completes
the proof of Theorem \ref{th1.6}.

\section{\bf New proof of sharp Michael-Simon inequality}
\noindent In this section, we will apply the smooth convergence
result in Theorem \ref{th1.6} to give a new proof of the sharp
Michael-Simon inequality (\ref{1.2}) and find the necessary and
sufficient condition for the establishment of the equality. We know
that the inequality (\ref{1.2}) can be simplified to (\ref{1.3}), so
the monotonicity of $\int_{M_{t}}f^{\frac{n}{n-1}}d\mu_{t}$ is the
key point.
\begin{lemma}
    Along flow (\ref{1.5}), we have
    \begin{align*}
        \frac{\partial}{\partial t}\int_{M_{t}}f^{\frac{n}{n-1}}d\mu_{t} \leq 0
    \end{align*}
\end{lemma}

\noindent{\it \bf Proof.}~~
\begin{align*}
    \frac{\partial}{\partial t}\int_{M_{t}}f^{\frac{n}{n-1}}d\mu_{t} &=\int_{M_{t}}\left(\frac{n}{n-1}f^{\frac{1}{n-1}}\frac{\partial f}{\partial t} + f^{\frac{n}{n-1}}\Phi H \right)d\mu_{t}\\
    &=\int_{M_{t}}f^{\frac{1}{n-1}}\left(\frac{n}{n-1}\frac{\partial f}{\partial r}\sqrt{1+r^{-2}|\nabla r|^{2}}+f H\right)\Phi d\mu_{t}\\
\end{align*}
and $\Phi=-\left(\frac{n}{n-1}\frac{\partial f}{\partial
r}\sqrt{1+r^{-2}|\nabla r|^{2}} +f H \right)$, thus
\begin{align}\label{5.1}
    \frac{\partial}{\partial t}\int_{M_{t}}f^{\frac{n}{n-1}}d\mu_{t} &=-\int_{M_{t}}f^{\frac{1}{n-1}}\left(\frac{n}{n-1}\frac{\partial f}{\partial r}\sqrt{1+r^{-2}|\nabla r|^{2}}
    +f H\right)^{2} \leq 0
\end{align}
\hfill${\square}$

The monotonicity of $\int_{M_{t}}f^{\frac{n}{n-1}}d\mu_{t}$ yields
\begin{align*}
    \int_{M_{0}}f^{\frac{n}{n-1}}d\mu \geq \int_{M_{t}}f^{\frac{n}{n-1}}d\mu_{t} \geq \int_{M_{\infty}}f^{\frac{n}{n-1}}d\mu_{\infty}
\end{align*}
and from the smooth convergence of the flow (\ref{1.5}),
$M_{\infty}=B_{r_{\infty}}$, then
\begin{align*}
    H = nr_{\infty}^{-1};\qquad \nabla r_{\infty} =0
\end{align*}
and
\begin{align*}
    \nabla f(r_{\infty}) =\frac{\partial f}{\partial r_{\infty}} \nabla r_{\infty}= 0
\end{align*}
Therefore, $f(r_{\infty})$ is constant and
\begin{align}\label{5.2}
    \int_{M_{\infty}}f^{\frac{n}{n-1}}d\mu_{\infty} =\int_{B_{r_{\infty}}}f^{\frac{n}{n-1}}d\mu_{\infty}= f^{\frac{n}{n-1}}(r_{\infty})r_{\infty}^{n}|B^{n}|
\end{align}
The same, we have $\frac{\partial}{\partial
t}\int_{M_{t}}f^{\frac{n}{n-1}}d\mu_{t}=0$ when $M_{t} =
B_{r_{\infty}}$, then
\begin{align}\label{5.3}
    \frac{n}{n-1}\frac{\partial f}{\partial r}\sqrt{1+r^{-2}|\nabla r|^{2}}
     +f H = 0
\end{align}
i.e.
\begin{align*}
    nf(r_{\infty})r_{\infty}^{-1}+\frac{n}{n-1}\frac{\partial f}{\partial r}(r_{\infty}) = 0
\end{align*}
Solving the above ODE, we infer that
\begin{align}\label{5.4}
    f(r_{\infty})= r_{\infty}^{-(n-1)}
\end{align}
and adding (\ref{5.4}) into (\ref{5.2}), we get
\begin{align*}
    \int_{M_{\infty}}f^{\frac{n}{n-1}}d\mu_{\infty}= |B^{n}|
\end{align*}
hence
\begin{align*}
    \int_{M} f^{\frac{n}{n-1}} d\mu \geq |B^{n}|
\end{align*}
which means that the starshaped hypersurface satisfies the
inequality (\ref{1.2}).

It is obvious that equality holds for the sphere, we just need to
prove the converse. Suppose that the smooth starshaped hypersurface
$M_{t}$ makes the equality hold
\begin{align}\label{5.5}
    \int_{M_{t}}f^{\frac{n}{n-1}}d\mu= f^{\frac{n}{n-1}}(r_{\infty})r_{\infty}^{n}|B^{n}|=|B^{n}|
\end{align}
Then, along the flow (\ref{1.5}), the integral
$\int_{M_{t}}f^{\frac{n}{n-1}}d\mu_{t}$ remains to be a constant and
there is (\ref{5.3}) that
\begin{align*}
  \frac{n}{n-1}\frac{\partial f}{\partial r}\sqrt{1+r^{-2}|\nabla r|^{2}} +f H = 0
\end{align*}
Thereby, we find
\begin{align*}
    \frac{\partial r}{\partial t} = -\sqrt{1+r^{-2}|\nabla r|^{2}}\left(\frac{n}{n-1}\frac{\partial f}{\partial r}\sqrt{1+r^{-2}|\nabla r|^{2}}
    +f H \right) =0
\end{align*}
and
\begin{align*}
    \frac{\partial}{\partial t}\left(f^{\frac{n}{n-1}}r^{n}\right) = f^{\frac{1}{n-1}}\left(\frac{n}{n-1}\frac{\partial f}{\partial r}r+nf\right)r^{n-1}\frac{\partial r}{\partial t} =0
\end{align*}
Thus
\begin{align*}
    r^{n}f^{\frac{n}{n-1}}(r) =r_{\infty}^{n}f^{\frac{n}{n-1}}(r_{\infty})
\end{align*}
The equation (\ref{5.5}) is equivalent to
\begin{align}\label{5.6}
    \int_{M_{t}}f^{\frac{n}{n-1}}d\mu= f^{\frac{n}{n-1}}(r_{\infty})r_{\infty}^{n}|B^{n}|=f^{\frac{n}{n-1}}(r)r^{n}|B^{n}|
\end{align}
and we can deduce that
\begin{align}\label{5.7}
    r^{n} = \frac{\int_{M_{t}}f^{\frac{n}{n-1}}d\mu_{t}}{f^{\frac{n}{n-1}}(r)|B^{n}|}
\end{align}
Differentiating (\ref{5.7}) and combining with (\ref{5.6}), we get
\begin{align*}
    \nabla(r^{n}) = \nabla \left( \frac{\int_{M_{t}}f^{\frac{n}{n-1}}d\mu_{t}}{f^{\frac{n}{n-1}}(r)|B^{n}|} \right)
    = - \frac{n}{n-1}\frac{\nabla f}{f}r^{n}
\end{align*}
Thus
\begin{align*}
    \nabla(r^{n})+\frac{n}{n-1}\frac{\nabla f}{f}r^{n} =\nabla r \left(f+\frac{n}{n-1}\frac{\partial f}{\partial r} r\right) = 0
\end{align*}
It can be inferred that either $\nabla r = 0$ or
$f+\frac{n}{n-1}\frac{\partial f}{\partial r} r=0$.

Case1: If $\nabla r = 0$, $r$ is constant, then $M_{t}=B_{r}$ and
$\nabla f= \frac{\partial f}{\partial r}\nabla r = 0$ i.e. $f$ is
constant.

Case2: If $f+\frac{1}{n-1}\frac{\partial f}{\partial r} r=0$, then
\begin{align}\label{5.8}
    f(r)= r^{-(n-1)}
\end{align}
Taking (\ref{5.8}) into (\ref{5.3}), we have
\begin{align*}
    \left( -\frac{1}{n-1}rH + \frac{n}{n-1}\sqrt{1+r^{-2}|\nabla r|^{2}}\right)\frac{\partial f}{\partial r} = 0
\end{align*}
i.e.
\begin{align*}
    -\frac{1}{n-1}rH + \frac{n}{n-1}\sqrt{1+r^{-2}|\nabla r|^{2}}=0
\end{align*}
Furthermore, we can get
\begin{align}\label{5.9}
    H= \frac{n\sqrt{r^{2}+|\nabla r|^{2}}}{r^{2}}
\end{align}
Owing to
\begin{align}\label{5.10}
    H=g^{ij}h_{ij}= \frac{n(r^{2}+2|\nabla r|^{2}-r\Delta r)}{(r^{2}+2|\nabla r|^{2})^{\frac{3}{2}}}
\end{align}
Combining (\ref{5.9}) and (\ref{5.10}), we have
\begin{align*}
    r \Delta r = -r^{-2}|\nabla r|^{4}
\end{align*}
and the equation (\ref{2.2}) implies that 
\begin{align*}
    h_{ij} = \frac{e_{ij}}{\sqrt{r^{2}+|\nabla r|^{2}}}\left( r+r^{-1}|\nabla r|^{2}\right)^{2}
\end{align*}
Thus
\begin{align*}
    |A|^{2}=g^{ij}g^{kl}h_{ik}h_{il} = \frac{n^{2}}{(r^{2}+|\nabla r|^{2})^{3}} (r+ r^{-1}|\nabla r|^{2})^{4}
\end{align*}
Note that
\begin{align*}
    \frac{|A|^{2}}{H^{2}} = 1
\end{align*}
Therefore, $M_{t}$ is sphere and it can be inferred that $f$ is
constant.
\section{\bf Preserving of static convexity}
\noindent In this section, we will show that static convexity is
preserved along the flow (\ref{1.7}). The main tool we use in our
proof is the tensor maximum principle shown below, which was proved
by Andrews.
\begin{theorem}(\cite{And07})\label{th6.1}
    Let $S_{i j}$ be a smooth time-varying symmetric tensor field on a compact manifold $M$ satisfying
    \begin{align*}
        \frac{\partial}{\partial t} S_{i j}=a^{k l} \nabla_{k} \nabla_{l} S_{i j}+b^{k} \nabla_{k} S_{i j}+N_{i j}
    \end{align*}
    where $a^{k l}$ and $b$ are smooth, $\nabla$ is a (possibly time-dependent) smooth symmetric connection, and $a^{k l}$ is positive definite everywhere. Suppose that
    \begin{align}\label{6.1}
        N_{i j} v^{i} v^{j}+\sup _{\Lambda_{k}^{p}} 2 a^{k l}\left(2 \Lambda_{k}^{p} \nabla_{l} S_{i p} v^{i}-\Lambda_{k}^{p} \Lambda_{l}^{q} S_{p q}\right) \geq 0
    \end{align}
    whenever $S_{i j} \geq 0$ and $S_{i j} v^{j}=0$, where the supremum is taken over all $n \times n$ matrix $\Lambda_{k}^{p} .$ If $S_{i j} \geq 0$ everywhere on $M$ at $t=0$ and on $\partial M$ for $0 \leq t \leq T$, then $S_{i j} \geq 0$ holds on $M$ for $0 \leq t \leq T$.
\end{theorem}
Now, we state the main conclusions of this section.
\begin{theorem}\label{th6.2}
   If the initial hypersurface $M_{0}$ is static convex, then along the flow (\ref{1.7}), the evolving hypersurfaces $M_{t}$ is static convex for $t>0$.
\end{theorem}

\noindent{\it \bf Proof.} Let $S_{ij}=h_{ij}-h^{-1}g_{ij}$, then
static convexity is equivalent to the positive of the tensor
$S_{ij}$. By (\ref{3.7}), we have
\begin{align*}
    \frac{\partial}{\partial t}S_{ij}=&h\dot{F}^{kl}\nabla_{k}\nabla_{l}S_{ij}+h^{-2}\dot{F}^{kl}\nabla_{k}h\nabla_{l}hg_{ij}+h\ddot{F}^{kl,pq}\nabla_{i}h_{kl}\nabla_{j}h_{pq}+\nabla_{i}h\nabla_{j}F+\nabla_{j}h\nabla_{j}F \\
    &+(1-3hF)(S^{2})_{ij}+\left(h\dot{F}^{kl}(h^{2})_{kl}-3F\right)S_{ij} +2\dot{F}^{kl}(h^{2})_{kl}g_{ij}-2h^{-1}Fg_{ij}
\end{align*}
From Theorem \ref{th6.1}, we need to prove the inequality
(\ref{6.1}) whenever $S_{ij}\geq 0$ and $S_{ij}v^{i}=0$ ($v$ is a
null vector of $S_{ij}$). Suppose that at $(x_{\tau},\tau)$,
$S_{ij}$ has a null vector $v$ and the principal curvatures is
strictly monotonically increasing, i.e.
$\kappa_{1}<\kappa_2<\cdots<\kappa_{n}$. $S_{ij}v^{i}=0$ implys that
$v=e_{1}$ and $S_{11}=\kappa_{1}-h^{-1}=0$ at $(x_{\tau},\tau)$.
Since $S_{ij}$ and $(S^{2})_{ij}$ satisfy the null vector condition,
it remains to show that
\begin{align*}
    Q_{1}:=&2h^{-2}\dot{F}^{kl}\nabla_{k}h\nabla_{l}hg_{ij}+h\ddot{F}^{kl,pq}\nabla_{i}h_{kl}\nabla_{j}h_{pq}+\nabla_{i}h\nabla_{j}F+\nabla_{j}h\nabla_{j}F\\
    &+2\dot{F}^{kl}(h^{2})_{kl}g_{ij}-2h^{-1}Fg_{ij} +2h \sup_{\varLambda }\dot{F}^{kl}\left(2\varLambda^{p}_{q}\nabla_{l}S_{1p}-\varLambda^{p}_{k}\varLambda^{q}_{l}S_{pq}\right)\\
    \geq & 0
\end{align*}
at $\left(x_{\tau},\tau\right)$ for all matrix
$\left(\Lambda_{k}^{p}\right)$. We know that $S_{11}=0$ and
$\nabla_{k} S_{11}=0$ at $\left(x_{\tau}, \tau\right)$, then
\begin{align*}
\dot{F}^{k l}\left(2 \Lambda_{k}^{p} \nabla_{l} S_{1 p}-\Lambda_{k}^{p} \Lambda_{l}^{q} S_{p q}\right) &=\sum_{k=1}^{n} \sum_{p=2}^{n} \dot{F}^{k}\left(2 \Lambda_{k}^{p} \nabla_{k} S_{1 p}-\left(\Lambda_{k}^{p}\right)^{2} S_{p p}\right) \\
&=\sum_{k=1}^{n} \sum_{p=2}^{n}
\dot{F}^{k}\left(\frac{\left(\nabla_{k} S_{1 p}\right)^{2}}{S_{p
p}}-\left(\Lambda_{k}^{p}-\frac{\nabla_{k} S_{1 p}}{S_{p
p}}\right)^{2} S_{p p}\right)
\end{align*}
It follows that the supremum is obtained by choosing
$\Lambda_{k}^{p}=\frac{\nabla_{k} S_{1 p}}{S_{p p}}$ for $p \geq 2$
, $ k \geq 1$ and $\Lambda_{k}^{1}=0$ for all $k$. Thus, $Q_{1}$
becomes:
\begin{align*}
    Q_{1}:=&2h^{-2}\dot{F}^{kl}\nabla_{k}h\nabla_{l}hg_{ij}+h\ddot{F}^{kl,pq}\nabla_{i}h_{kl}\nabla_{j}h_{pq}+\nabla_{i}h\nabla_{j}F+\nabla_{j}h\nabla_{j}F\\
    &+2\dot{F}^{kl}(h^{2})_{kl}g_{ij}-2h^{-1}Fg_{ij}+2h\sum_{k=1}^{n} \sum_{p=2}^{n} \dot{F}^{k}\frac{\left(\nabla_{k} S_{1 p}\right)^{2}}{S_{p p}}
\end{align*}
By Lemma \ref{le2.6}, $F(\kappa)$ is inverse concave. Through direct
calculation, such as in \cite{AW18}, we have
\begin{align*}
   &h\ddot{F}^{kl,pq}\nabla_{1}h_{kl}\nabla_{1}h_{pq}+2h\sum_{k=1}^{n} \sum_{p=2}^{n} \dot{F}^{k}\frac{\left(\nabla_{k} S_{1 p}\right)^{2}}{S_{p p}}  \\
   &\geq 2hF^{-1}|\nabla_{1}F|^{2}+2\sum_{k>1,p>1}\dot{F}^{k}\left(\frac{1}{\kappa_{p}-h^{-1}}-\frac{1}{\kappa_{p}}\right) \left(\nabla_{1}h^{2}_{kp}\right)\\
   &\geq2hF^{-1}|\nabla_{1}F|^{2}+2\sum_{k>1}\dot{F}^{k}\frac{h^{-1}}{(\kappa_{k}-h^{-1})\kappa_{k}} \left(\nabla_{1}h^{2}_{kk}\right)
\end{align*}
In addition, using the Cauchy-Schwarz inequality, we have
\begin{align*}
    \sum_{k=2}^{n} \frac{\dot{F}^{k}}{\kappa_{k}\left(\kappa_{k}-h^{-1}\right)}\left(\nabla_{1} h_{k k}\right)^{2} \cdot \sum_{k=2}^{n} \dot{F}^{k} \kappa_{k}\left(\kappa_{k}-h^{-1}\right) \geq\left(\sum_{k=2}^{n} \dot{F}^{k}\left|\nabla_{1} h_{k k}\right|\right)^{2} \geq\left|\nabla_{1} F\right|^{2}
\end{align*}
where we use $\nabla_{1} S_{11}=\nabla_{1}
\kappa_{1}=\nabla_{1}h^{-1}=0$ and the inequalities (\ref{2.21})
that $\sum_{k=1}^{n} \dot{F}^{k} \kappa_{k}^{2}-F \geq F^{2}-F>0$.
Since $h^{-1}<\kappa_{2}<\cdots<\kappa_{n}$, we get
\begin{align*}
    \sum_{k=2}^{n} \frac{\dot{F}^{k}}{\kappa_{k}\left(\kappa_{k}-h^{-1}\right)}\left(\nabla_{1} h_{k k}\right)^{2} \geq \frac{\left|\nabla_{1} F\right|^{2}}{\sum_{k=1}^{n} \dot{F}^{k} \kappa_{k}^{2}-h^{-1}F}  \geq 0
\end{align*}
Thereby
\begin{align*}
    Q_{1}\geq& 2hF^{-1}|\nabla_{1}F|^{2}+2h\sum_{k>1}^{n}\dot{F}^{k}\frac{h^{-1}}{(\kappa_{k}-h^{-1})\kappa_{k}}(\nabla_{1}h_{kk})^{2}\\
    &+2h^{-2}\dot{F}^{kl}\nabla_{k}h\nabla_{l}h+2\dot{F}^{kl}(h^{2})_{kl}-2h^{-1}F\\
    \geq& 2hF^{-1}|\nabla_{1}F|^{2}+2\frac{|\nabla_{1}F|^{2}}{\sum_{k=1}^{n}\dot{F}^{k}\kappa_{k}^{2}-h^{-1}F}+2\dot{F}^{kl}(h^{2})_{kl}-2h^{-1}F\\
    &+2h^{-2}\dot{F}^{kl}\langle X, x_{i}\rangle \langle X, x_{j}\rangle h^{i}_{k}h^{j}_{l}\\
    \geq& 2h^{-2}F^{2}+2F^{2}-2h^{-1}F\\
    \geq& 2h^{-1}F\left(h^{-1}F-1\right)+2F^{2} \geq 0
\end{align*}

\hfill${\square}$

\noindent {\bf Proof of Theorem \ref{th1.9}}: We proved that static
convexity is preserved along the flow (\ref{1.7}), then $M_{t}$ is
strictly convex for all $t>0$ and can be represented  $M_{t}$ in
terms of the support function $h$. First of all, directly applying
the maximum principle we can obtain $C^{0}$ estimate of $h$, which
is equivalent to $C^{0}$ estimate of the solution $M_{t}$. $C^{1}$
estimate and $C^{2}$ estimate follow from the same steps in
\cite{GL2021}. Secondly, since the initial hypersurface $M_{0}$ is
static convex, which implies that $M_{0}$ is strictly convex, then
one can obtain that flow (\ref{1.7}) is uniformly parabolic. Through
Krylov’s theory, we get $C^{2,\alpha}$ estimate. The longtime
existence of the flow (\ref{1.7}) can be obtained from standard
theory for parabolic equations. Finally, as in \cite{GL2018}, By
estimating the gradient of the radial function of $M_{t}$, we can
derive that $M_{t}$ is exponentially converging to a sphere.
\section{\bf Sharp Michael-Simon inequality for \texorpdfstring{$k$} th mean curvature}
\noindent In this section, we use the convergence result of the flow
(\ref{1.7}) to prove the inequality (\ref{1.10}) for static convex
hypersurface. First of all, by scaling, we may assume that
\begin{align*}
   \int_{M} \sqrt{\sigma_{k}^{2}f^{2}+\sigma_{k-1}^{2}|\nabla^{M} f|^{2}} +\int_{\partial M} \sigma_{k-1}f = n\int_{M}\sigma_{k-1}f^{\frac{n-k+1}{n-k}}
\end{align*}
This normalization ensures that we can find a function $\vartheta: M
\rightarrow \mathbb{R}$ which solves the PDE
\begin{align*}
    \operatorname{div}_{M}\left(\sigma_{k-1}f \nabla^{M} \vartheta \right)=n \sigma_{k-1}f^{\frac{n-k+1}{n-k}}-\sqrt{\sigma_{k}^{2}f^{2}+\sigma_{k-1}^{2}|\nabla^{M} f|^{2}}
\end{align*}
on $M$ with Neumann boundary condition $\left\langle\nabla^{M}
\vartheta , \vec{n}\right\rangle=1$ on $\partial M$. Here, $\vec{n}$
denotes the co-normal to $M$. Note that $\vartheta$ is of class
$C^{2, \beta}$ for each $0<\beta<1$ by standard elliptic regularity
theory. Now we only need to prove the following inequality
\begin{align}\label{7.1}
    \int_{M}\sigma_{k-1}f^{\frac{n-k+1}{n-k}} d\mu \geq y_{k}\circ z_{k-1}^{-1}(V_{k-1}(\Omega))
\end{align}

Secondly, we prove that the flow (\ref{1.7}) preserves the $(k-1)$th
quermassintergral $V_{k-1}(\Omega_{t})$ where $\Omega_{t}$ is the
domain enclosed by $M_{t}$, as following
\begin{align*}
    \frac{\partial}{\partial t}V_{k-1}(\Omega_{t})=(n+2-k)\int_{M_{t}}\left(1-h\frac{E_{k}}{E_{k-1}}\right)E_{k-1}d\mu_{t}=\int_{M_{t}}\left(E_{k-1}-hE_{k}\right)d\mu_{t}=0
\end{align*}
where the last equality we use Minkowski identity (\ref{2.18}).

Last but not least, we need to establish the monotonicity of
$\int_{M_{t}} \sigma_{k-1}(\kappa)f^{\frac{n-k+1}{n-k}} d\mu_{t}$
along the flow (\ref{1.7}). For computational convenience, using
$g=f^{\frac{n-k+1}{n-k}}$
\begin{align*}
   \frac{\partial}{\partial t}\int_{M_{t}}\sigma_{k-1}g d\mu_{t}=&\int_{M_{t}}\left(g \frac{\partial \sigma_{k-1}}{\partial t}+\sigma_{k-1}\frac{\partial g}{\partial f}\frac{\partial f}{\partial h}\frac{\partial h}{\partial t} + \sigma_{1}\sigma_{k-1}g\Phi  \right)d\mu_{t} \\
   =&\int_{M_{t}} g\binom{n}{k-1} \dot{E}^{ij}_{k-1}\left(-\nabla_{j}\nabla_{i}\Phi - \Phi(h^{2})_{ij}\right)d\mu_{t}\\
   & + \int_{M_{t}}\left(\sigma_{k-1}\frac{\partial g}{\partial f}\frac{\partial f}{\partial h}+g\sigma_{k-1}\sigma_{1}\right)\Phi d\mu_{t}\\
   =&\int_{M_{t}}-\binom{n}{k-1}\left(\dot{E}^{ij}_{k-1}\nabla_{j}\nabla_{i}g\right)\Phi d\mu_{t}+\int_{M_{t}}\frac{\partial g}{\partial f}\frac{\partial f}{\partial h}\sigma_{k-1}\Phi d\mu_{t}\\
   &+\int_{M_{t}}kg\sigma_{k}\Phi d\mu_{t}\\
   =&\int_{M_{t}} -\binom{n}{k}(k-1)\left(\frac{\partial^{2}g}{\partial f^{2}}\left(\frac{\partial f}{\partial h}\right)^{2}+\frac{\partial g}{\partial f}\frac{\partial^{2}f}{\partial h^{2}}\right)\dot{E}^{ij}_{k-1}\nabla_{i}h\nabla_{j}h\Phi d\mu_{t}\\
   &-\int_{M_{t}}\binom{n}{k-1}\frac{\partial g}{\partial f}\frac{\partial f}{\partial h}\dot{E}^{ij}_{k-1}\nabla_{i}\nabla_{j}h\Phi d\mu_{t}\\
   &+ \int_{M_{t}}k\binom{n}{k}gE_{k}\Phi d\mu_{t}
   +\int_{M_{t}}\binom{n}{k-1}\frac{\partial g}{\partial f}\frac{\partial f}{\partial h}E_{k-1}\Phi d\mu_{t}
\end{align*}
where the third equality we use integration by parts. Combining
(\ref{2.13}) and (\ref{2.14}), we have
\begin{align}
    \frac{\partial}{\partial t}\int_{M_{t}}\sigma_{k-1}g d\mu_{t}=&\int_{M_{t}} -\binom{n}{k}(k-1)\left(\frac{\partial^{2}g}{\partial f^{2}}\left(\frac{\partial f}{\partial h}\right)^{2}+\frac{\partial g}{\partial f}\frac{\partial^{2}f}{\partial h^{2}}\right)\dot{E}^{ij}_{k-1}(h^{2})_{ij}\Phi d\mu_{t}\notag\\
    &-\int_{M_{t}}\binom{n}{k-1}\frac{\partial g}{\partial f}\frac{\partial f}{\partial h}\dot{E}^{ij}_{k-1}\left(h_{ij}-h(h^{2})_{ij}\right)\Phi d\mu_{t} \notag\\
    &+ \int_{M_{t}}k\binom{n}{k}gE_{k}\Phi d\mu_{t}
    +\int_{M_{t}}\binom{n}{k-1}\frac{\partial g}{\partial f}\frac{\partial f}{\partial h}E_{k-1}\Phi d\mu_{t}
    \notag\\
    =&\int_{M_{t}} -\binom{n}{k}(k-1)\left(\frac{\partial^{2}g}{\partial f^{2}}\left(\frac{\partial f}{\partial h}\right)^{2}+\frac{\partial g}{\partial f}\frac{\partial^{2}f}{\partial h^{2}}\right)\dot{E}^{ij}_{k-1}(h^{2})_{ij}\Phi d\mu_{t} \notag\\
    &+\int_{M_{t}}\binom{n}{k-1}\frac{\partial g}{\partial f}\frac{\partial f}{\partial h}h\dot{E}^{ij}_{k-1}(h^{2})_{ij}\Phi d\mu_{t}\notag\\
    &+ \int_{M_{t}}k\binom{n}{k}gE_{k}\Phi d\mu_{t}
    +\int_{M_{t}}\binom{n}{k-1}\frac{\partial g}{\partial f}\frac{\partial f}{\partial h}E_{k-1}\Phi d\mu_{t} \notag\\
    &-\int_{M_{t}}\binom{n}{k-1}(k-1)\frac{\partial g}{\partial f}\frac{\partial f}{\partial h}
    E_{k-1}\Phi d\mu_{t} \label{7.2}
\end{align}
where we use (\ref{2.16}) in the last equality. Adding $\Phi =
1-h\frac{E_{k}}{E_{k-1}}$ in (\ref{7.2}), we get
\begin{align*}
    \frac{\partial}{\partial t} & \int_{M_{t}} \sigma_{k-1}g d\mu_{t}\\
    =&\int_{M_{t}}-\binom{n}{k}(k-1)\left(\frac{\partial^{2}g}{\partial f^{2}}\left(\frac{\partial f}{\partial h}\right)^{2}+\frac{\partial g}{\partial f}\frac{\partial^{2}f}{\partial h^{2}}\right)\left(nE_{1}-(n-k+1)\frac{E_{k}}{E_{k-1}}\right) (E_{k-1}-hE_{k}) d\mu_{t}\\
    &-\int_{M_{t}}\binom{n}{k-1}(k-2)\frac{\partial g}{\partial f}\frac{\partial f}{\partial h}(E_{k-1}-hE_{k}) d\mu_{t} +
    \int_{M_{t}}k\binom{n}{k}g\frac{E_{k}}{E_{k-1}}(E_{k-1}-hE_{k})
    d\mu_{t} \\
    &+\int_{M_{t}} n\binom{n}{k-1}\frac{\partial g}{\partial f}\frac{\partial f}{\partial h}hE_{1}(E_{k-1}-hE_{k}) d\mu_{t}\\
    &-\int_{M_{t}}(n-k+1)\binom{n}{k-1}\frac{\partial g}{\partial f}\frac{\partial f}{\partial h}\frac{E_{k}}{E_{k-1}}(E_{k-1}-hE_{k}) d\mu_{t}
\end{align*}
where we use (\ref{2.16}). From Assumption \ref{as1.10}, we can get
the following inequalities
\begin{align*}
\frac{\partial f}{\partial h} \geq 0;\qquad
\frac{\partial^{2}g}{\partial f^{2}}\left(\frac{\partial f}{\partial
h}\right)^{2}+\frac{\partial g}{\partial
f}\frac{\partial^{2}f}{\partial h^{2}} \leq 0
 \end{align*}
and using (\ref{1.9}) and (\ref{2.19}), we also have   
\begin{align*}
    h_{ij} \geq h^{-1}g_{ij} \left(i.e. E_{1}\geq h^{-1}\right); \quad
   (E_{k-1}-hE_{k})=\frac{1}{k-1}\dot{E}_{k-1}^{ij}(g_{ij}-hh_{ij})\leq 0;\quad E_{k}\leq E_{1}E_{k-1}
\end{align*}
Thus, combining the above several inequalities, we get
\begin{align}
   \frac{\partial}{\partial t}\int_{M_{t}}\sigma_{k-1}g d\mu_{t}\leq& \int_{M_{t}} -\binom{n}{k-1}(k-2)\frac{\partial g}{\partial f}\frac{\partial f}{\partial h}(E_{k-1}-hE_{k})d\mu_{t} \notag\\
   &-\int_{M_{t}}(n-k+1)\binom{n}{k-1}\frac{\partial g}{\partial f}\frac{\partial f}{\partial h}hE_{1}(E_{k-1}-hE_{k})d\mu_{t} \notag\\
   &+ \int_{M_{t}}n\binom{n}{k-1}\frac{\partial g}{\partial f}\frac{\partial f}{\partial h}hE_{1}(E_{k-1}-hE_{k})d\mu_{t} \notag\\
   \leq &\int_{M_{t}} -\binom{n}{k-1}(k-2)\frac{\partial g}{\partial f}\frac{\partial f}{\partial h}(E_{k-1}-hE_{k})d\mu_{t}  \notag\\
   &+\int_{M_{t}}(k-1)\binom{n}{k-1}\frac{\partial g}{\partial f}\frac{\partial f}{\partial h}hE_{1}(E_{k-1}-hE_{k})d\mu_{t} \notag\\
   \leq & \int_{M_{t}}\binom{n}{k-1}\frac{\partial g}{\partial f}\frac{\partial f}{\partial h}(E_{k-1}-hE_{k})d\mu_{t} \leq 0  \label{7.3}
\end{align}
Therefore
\begin{align*}
   \int_{M_{t}} \sigma_{k-1}g d\mu_{t} \geq \int_{M_{\infty}}\sigma_{k-1}g d\mu_{\infty}=\int_{B_{R}}\sigma_{k-1}g d\mu
\end{align*}
The equality in the above formula is obtained from the convergence
result of the flow (\ref{1.7}) and $B_{R}=\partial B^{n+1}_{R}$.
Also $\nabla f=\frac{\partial f}{\partial h}\nabla h=0$ on $B_{R}$
i.e. $f$ is constant on $B_{R}$.

We already know the flow (\ref{1.7}) preserved $(k-1)$th
quermassintergral, i.e.
\begin{align*}
    V_{k-1}(\Omega_{0})=V_{k-1}(\Omega_{t})=V_{k-1}(B^{n+1}_{R})
\end{align*}
and
\begin{align*}
   \int_{B_{R}} \sigma_{k-1}g d\mu &=\binom{n}{k-1}f^{\frac{n-k+1}{n-k}}\int_{B_{R}} E_{k-1}(\kappa) d\mu=\binom{n}{k-1}f^{\frac{n-k+1}{n-k}}V_{k}(B^{n+1}_{R})\\
   &=y_{k}\circ z^{-1}_{k-1}(V_{k-1}(B^{n+1}_{R}))=y_{k}\circ z^{-1}_{k-1}(V_{k-1}(\Omega_{t}))=y_{k}\circ z^{-1}_{k-1}(V_{k-1}(\Omega_{0}))
\end{align*}
Then
\begin{align}\label{7.4}
    \int_{M_{t}} \sigma_{k-1}g d\mu_{t} \geq y_{k}\circ z^{-1}_{k-1}(V_{k-1}(\Omega_{t}))
\end{align}

Now we just need to show that $M$ is a sphere and $f$ is constant
when the equality holds in (\ref{1.10}). If the smooth static convex
hypersurface $M_{t}$ attains the equality
\begin{align*}
    \int_{M_{t}}\sigma_{k-1}f^{\frac{n-k+1}{n-k}} d\mu_{t} = y_{k}\circ z^{-1}_{k-1}(V_{k-1}(\Omega_{t}))
\end{align*}
Then the equality in (\ref{7.3}) holds, which implies that the
principal curvature of $M_{t}$ satisfies
$\kappa_{1}=\cdots=\kappa_{n}=h^{-1}$ i.e. $M$ is the sphere. In the
same way, we can infer that $f$ is constant. This completes the
proof of Theorem \ref{th1.11}.

Without loss of generality, if M is a sphere $B_{R}$ and we take
$f=R^{-(n-k)}$, then
\begin{align*}
    \int_{B_{R}}\sigma_{k-1}f^{\frac{n-k+1}{n-k}} d\mu = |B^{n}|
\end{align*}
Therefore, we obtain the inequality (\ref{1.4}) for the static
convex hypersurface.

\section*{\bf Acknowledgments}
 The authors would like to thank
Professors X. Yang,  X. Zhang and X. Chao for their guidance and
help! This work is supported by NNSF of China(no.11871275).


\end{document}